\documentclass[a4paper,10pt]{article}

\usepackage{amssymb,amsfonts}
\usepackage{amsmath}

\usepackage{graphicx}
\begin{document}
\title{Some remarks on polar sets to sums of squares}
\author{Tove Dahn (Lund University)}

\maketitle

\section{Introduction}
The problem we consider, is to determine the character of the movements, that preserve hypoellipticity
for symbols represented as sums of squares. The main result is that a necessary condition for 
hypoellipticity, is that the polar is not a spiral domain. A sufficient condition is that
the polar is algebraic.

$\mathcal{G}$ denotes groups of transformations (movements) defined according to (\cite{Lie91})
by $X(f)=\xi \frac{\delta f}{\delta x} + \eta \frac{\delta f}{\delta y}$, where $f \in \mathcal{D}_{L^{1}}$ (\cite{Schwartz66}).
We assume $\mathcal{G}$ is very regular, in the sense that $\exists V \in \mathcal{G}$,
such that $V f$ is analytic, that is $X_{V}(f)=0$, on a domain of $x,y$. The movements  $U,V \in
\mathcal{G}$ are assumed to be related through $\frac{d U}{d V}=\alpha$, where $\alpha$
is polynomial close to the boundary and that we have linear independence in the infinity, 
(preserves constant value (\cite{Cousin95}), collar point (\cite{Helffer})). The boundary is the points of $x,y$ where the movement
changes character $U \rightarrow V$. The movements are assumed dependent on one parameter and
such that when the movement is sequential, it preserves polynomials.

\subsection{The concept of movement}
 Given $V$ a movement
 such that $V^{-1} \overline{\delta}(V f)=X_{V}(f)=0$, $V \in \mathcal{G}$, then given $f$ 
 analytic, we have that $Vf$ is analytic,
 that is $Vf$ is holomorphic or constant. 
 
 \newtheorem{adj}{Lemma The associated transformation}[section]
 \begin{adj}
 Assume $U=U_{1}$, then we have $U (\frac{d f}{d x})(x,y)=(\frac{d f}{d x})(x+a,y+b)$, for constants $a,b$. 
 Further, $\frac{d}{d x}(Uf)=\frac{df}{d x}(x+ a,y+b)$, that is $U^{-1} \frac{d}{d x} Uf = \xi \frac{d f}{d x}$ and so on. 
 In particular $X_{V}(f)=V^{-1} \overline{\delta} V f=g$ or $V g=\overline{\delta} V f$, why
 if $V \rightarrow I$, we assume $\overline{\delta} f=g$. Given $F$ a Hamiltonian, 
 we have
 $X(F)=0$ iff $\eta / \xi \sim - Y / X$. Note that if $(X,Y) \bot (\xi,\eta)$, 
 that is $\xi X + \eta Y=0$, we have $- \eta \frac{\delta f}{\delta x} + \xi \frac{\delta f}{\delta y}=0$.
 \end{adj}
 
  Assume $(X,Y) \rightarrow (-Y,X)$ projective, that is bijective and maps zero lines on zero lines (\cite{Lie91}).
 Let $X^{\diamondsuit}(f)=-\eta \frac{\delta f}{\delta x} + \xi \frac{\delta f}{\delta y}$. 
 Note that $(-\eta + i \xi)=i(\xi + i \eta)$. Thus, $X \rightarrow X^{\diamondsuit}$ preserves 
 analyticity, given $\xi,\eta$ such that $\eta / \xi \sim - Y / X$.

  \newtheorem{rorelse}[adj]{Definition Harmonic generator}
 \begin{rorelse}
 Assume that $\Omega$ is defined by $\mathcal{G}$, that is $(I)=\{ V f_{0} \quad V \in \mathcal{G} \}$
 and $\Omega=N(I)$, with $f_{0}$ holomorphic. For a fixed $V$, let $A=\{ V f_{0}=0 \}$ and $X_{V}(f_{0})=0$. 
 The condition $A \cap U \subset \Omega$, for a small set $U$, means that $(V + V_{1})f_{0}=0$, 
 where $V_{1}$ is assumed continuous, that is $\mbox{ reg }A$ corresponds to $V_{1} f_{0} \rightarrow 0$, 
 where $V_{1} \in \mathcal{G}$ and generates $U$, but is not necessarily analytic. In the same manner
 for $A^{\diamondsuit}$, where $A \rightarrow A^{\diamondsuit}$ is projective.
 \end{rorelse}

 \newtheorem{transversal}[adj]{Lemma The Hamiltonian defines a harmonic orthogonal}
 \begin{transversal}
 Given $Af=X \frac{\delta f}{\delta x} + Y \frac{\delta f}{\delta y}$, we have on 
 a contractible domain, where $f$ is Hamiltonian, that $A f \equiv 0$. We can define a 
 movement $U^{\diamondsuit}$ such that $U^{\diamondsuit -1}\frac{\delta}{\delta x}(U^{\diamondsuit} g)=-Y \frac{\delta g}{\delta x}$,
 $U^{\diamondsuit-1}\frac{\delta}{\delta y}(U^{\diamondsuit} g)=X \frac{\delta g}{\delta y}$. When $(X,Y)$ analytic 
 with the condition $\overline{\delta}(U^{\diamondsuit} g)=0$, then $U^{\diamondsuit}$ 
 is analytic over g. Assume $A(g)=h$ with $\delta h=0$, then we have that $h$ is exact 
 on a contractible domain, that is $U^{\diamondsuit}$ is harmonic over g and $U,U^{\diamondsuit}$ are related 
 by duality.
 Thus, given $f$ has analytic derivatives, we have existence of $g$ and $U^{\diamondsuit}$ harmonic over 
 $g$, with analytic derivatives $-Y,X$.
 \end{transversal}
 
Given $\frac{\xi}{\eta} \rightarrow -\frac{\eta}{\xi}$ projective, we can relate 
$U \rightarrow U^{\diamondsuit}$ to reflection. Given 
$\frac{\delta g}{\delta x}=-\eta$ and $\frac{\delta g}{\delta y}=\xi$, we have that 
$\{ f,g \}=X(f)=0$. If $\frac{\delta f}{\delta x}=\frac{\delta h}{\delta y}$ and 
$\frac{\delta f}{\delta y}=\frac{\delta h}{\delta x}$, for instance $h={}^{t} f$, 
we have that $\{ h,f \}=X^{\diamondsuit}(f)=0$. Regularity properties for $\xi,\eta$ are
important for representation of projection operators. When the coefficients have unbounded
sub level sets, the corresponding transformation is considered as a functional.

 \newtheorem{harmonic}[adj]{Lemma Harmonic transformations}
 \begin{harmonic}
 When $\delta U^{-1} \overline{\delta} Uf =\delta X(f)$, we have that $\delta X(f)=0$ iff
 $(\xi_{x} + \xi_{y}' ) f_{x} + (\eta_{x}' + \eta_{y}'') f_{y} + \eta' f_{xy} + \xi ' f_{yx} + \xi f_{xx} + \eta'' f_{yy}=0$,
 where $\xi '$ denotes $i \xi$. 
 Thus we have $\delta X(f)=0$, if for instance $\frac{d y}{d x} \sim - \frac{(\eta_{x}' + \eta_{y}'')}{(\xi_{x} + \xi_{y}')}$
 and $\frac{\eta'}{\xi'} \sim - \frac{f_{yx}}{f_{xy}} \sim - \frac{f_{xx}}{f_{yy}}$.
 \end{harmonic}
 
  Translation can be given as $y'=\varphi(x + k y)$
 $x_{1}=x - k t$,$y_{1}= y + t$. $X( f) \equiv - k \frac{\delta f}{\delta x} + \frac{\delta f}{\delta y}$
 and $A f \equiv \frac{\delta f}{\delta x} - \varphi(x + k y) \frac{\delta f}{\delta y}$. 
 Rotation
 can be given by $\frac{x y'  - y}{x + y y'}=f(x^{2} + y^{2})$ and $X(\frac{y}{x}) \equiv 1 + (\frac{y}{x})^{2}$
 or $X (\mbox{ arctg } \frac{y}{x})=1$ (\cite{Lie91}, Ch. 6, Ex. 3)

 Consider $U f(\zeta)=U F(\gamma)(\zeta) \sim F({}^{t} U \gamma)(\zeta)$.
 
 \newtheorem{pred}[adj]{Definition Accessible points}
 \begin{pred}
 Assume $\Sigma_{j}=\{ (x,y) \quad U F(x,y)=U_{j} F(x,y) \}$ and $\Omega_{j}=\{ \zeta \quad Uf(\zeta)=U_{j}f(\zeta) \}$.
 Given ${}^{t} U$ analytic in $\zeta$, there is a continuous mapping, $\Sigma_{j} \rightarrow \Omega_{j}$. 
 Given $X(f)=0$, we have $Uf$ analytic. Given ${}^{t} U \simeq U^{\bot}=0$, we assume $U=I$, 
 that is we have multivalentness. A point $\zeta \in \Omega$, is said to be
 accessible for $\mathcal{G}$, if we have a chain of sequential movements
  $U_{j} \gamma \rightarrow U_{r} \gamma \rightarrow \zeta$, where the last movement is 
  analytic. The remaining movements can be regarded as functionals. That is,  we assume 
  ${}^{t} U_{r}$ analytic over $\gamma$ in a neighbourhood of $\zeta$, 
  for the remaining ones we assume $\frac{d {}^{t} U}{d U_{1}}={}^{t} \alpha$ and so on.
\end{pred} 
 
 Assume that a domain is generated by $d \mu$ $\mu=\mu_{1} \times \mu_{r}$ and 
 $\mu^{\bot} \in \mathcal{G}$ with $\mu_{j} \neq \mu_{j}^{\bot}$. A normal model means 
 $\mu_{j} \mu_{j}^{\bot}=\mu_{j}^{\bot} \mu_{j}$,
 that is $\mu_{j} \mu_{k} = \mu_{k} \mu_{j}$ with $j \neq k$ and that the measure is finite. 
 For a spiral domain we have
 $\mu_{j}^{\bot}=\mu_{j}$, that is the condition on collar point is not satisfied. 
 We consider two completions of movement, $U \rightarrow U^{*} \rightarrow U^{\triangle}$ (Fourier dual)
 and $U \rightarrow \overline{U}^{\diamondsuit} \rightarrow U^{\triangledown}$ (harmonic conjugation) and finally,
 $U^{\blacktriangle}$ denotes $I-U$. The more general notation $U^{\bot}$ is relative context.
 For instance, relative the scalar product over $f,g \in \mathcal{D}_{L^{1}}$,
 $<f,\widehat{g}>$, we have  $<U f, \widehat{g}> \rightarrow <f, U^{\triangle} \widehat{g}>$
 and when both sides are $0$, we write ${}^{t} U \simeq U^{\bot}$.
 Given a parameter space, we can consider $U_{j}^{\blacktriangle}=0$, $j=1,\ldots,n$ as a foliation.
 The Stieltjes integral is then undetermined when $U_{j}^{\blacktriangle}=0$. Given a movement $U_{1}$ 
 such that $U_{1}^{\bot} \widehat{f}=0$ and $U_{1} \widehat{f}=0$ implies $U_{1}=I$, we have 
 ''projectivity``.

 \section{Involution}
 
 \subsection{Exact forms}

 Given $Uf$ surjective, 
 we have that ${}^{t} U (x,y)$ injective. Further, $\frac{d {}^{t} U}{d t}=\frac{\delta {}^{t} U}{\delta x} \frac{d x}{d t} + \frac{\delta {}^{t} U}{\delta y} \frac{d y}{d t}$.
 Assume that $\frac{\delta {}^{t} U}{\delta x}=\eta$,$\frac{\delta {}^{t} U}{\delta y}=-\xi$,
 why $\frac{d {}^{t} U}{d t}=\eta \frac{d x}{d t} - \xi \frac{d y}{d t}$, where we consider 
 $(x,y)(t)$, that is when $f$ is Hamiltonian we have that $\frac{d {}^{t} U}{d t}=-X(f)$. 
 When we assume $V^{-1} \overline{\delta}(V f)=X(f)$ with $V=U^{\diamondsuit}$, we can consider 
 $\int_{\beta} f d U^{\diamondsuit}=0$ that is vanishing flux. Assume $d U^{\diamondsuit}=p d \lambda$,
 for $d \lambda$ Lebesgue and $p$ a real polynomial.
 A boundary is given by $(x,y)$
 such that $f=0$ or $d U^{\diamondsuit}=0$. Given $f$ algebraic (cf. preserves constant value) 
 then the boundary is of measure zero. When $d U^{\diamondsuit}=0$ has positive measure, $U^{\diamondsuit}$
 preserves analyticity for $f$ (\cite{AhlforsSario60}).

 \newtheorem{dualitet}[adj]{Lemma Duality}
 \begin{dualitet}
 Assume $X(f)=g$ with $\overline{\delta}g=0$ and $X(f) \geq 0$.
 If we write $X(f) \sim <(\xi,\eta), d f>$ (more correctly $<(\xi,\eta),(-Y,X)>$), this defines 
 a linear functional on $d f$. 
 Thus given $f$ Hamiltonian, then $X(f)=0$ defines
 a duality $(\xi,\eta)$ to $(-Y,X)$. .
 
 \end{dualitet}

 For a global definition of the measures, it is sufficient to produce a continuous
 representation. Given a harmonic representation, using for instance
 $F(x,y) \rightarrow F(x,\frac{y}{x})$, where $y=y(x)$ continuous, we may have continuity
 in both $0,\infty$ simultaneously. Assume $\psi : (x,y) \rightarrow (x,\frac{y}{x})$.
 Assume $U_{\times}$ a sequential movement, $d U_{j} / d U_{1} = \alpha_{j}$ with $d U_{1}$ 
 harmonic. Given $\alpha_{j}^{\triangle} = e^{\phi^{\triangle}}$, with $\phi$ pseudo convex and 
 $\phi^{\triangle}$ a completion to harmonicity in the finite plane (\cite{Oka60}), 
 then we have that the movement changes character as $\mbox{ sgn }=$ $- / 0 / +$. 
 Assume $U_{\times}$ completed to $L^{1}$, then we have that $d U_{\times} f \in L^{1}$ 
 iff $\mid d U_{\times} f \mid \in L^{1}$. Let $\mid d U \mid \sim d U (\nu,\vartheta) $,
 where $\nu=y/x$ and $\vartheta=x/y$. 
 Note that given $(\nu,\vartheta)$ in a compact, we have that the same holds for $(x,y)$.
 If we only assume $\mid d U - d U^{\triangle} \mid \in L^{1}$,
 we have that $\nu \rightarrow \vartheta$ maps compact sets on compact sets.
 Assume $K$ a symmetric compact, such that $K \rightarrow 0$ through continuous deformation, 
 for instance $(x,y)$ such that $x^{2} + y^{2} < 2$ and $\mid x \mid < 1$ implies $\mid y/x \mid < 1$. 
 In the same manner if $\mid y \mid < 1$ we have that $\mid x / y \mid < 1$.
 Given $K$ is represented as closed (symmetric), we can assume $K$ a closed contour. 
 Thus, if any $K$ can be continuously deformed to $0$ using $\psi$, then $d U$ is a closed form. 
 Through the completion to harmonicity, the completed form can be seen as exact. 

 \newtheorem{exact}[adj]{Lemma Exactness for the completed form}
 \begin{exact}
 Given $\psi$ is proper, continuous with ${}^{t} \psi : \mathcal{D}_{L^{1}} \rightarrow \mathcal{D}_{L^{1}}$
 and such that every symmetric compact can be continuously deformed to a point, 
 the completed closed form can be seen as exact. 
 \end{exact}
 
   Assume $d U \sim \xi d x + \eta d y$, then we have that 
  $\int_{\Gamma} \xi  d x + \eta d y \sim \int_{(\Gamma)} (\eta_{x} - \xi_{y}) dx d y$
  thus for a harmonic movement, $(\Gamma)$ is contractible. Given $ d U = \alpha d U_{1}$, 
  where $\alpha \neq 0$ on $(\Gamma)$, the domain is contractible for $d U$.
 
\subsection{Evolute}
 Assume $\Sigma_{1}$ a given curve (involute) and $\Sigma_{2}$ its evolute (polar set according to Lie),
 then we have that $(T \Sigma_{1})^{\bot} \simeq T \Sigma_{2}$ (envelop of normals). 
 The determination of evolutes can be reduced to determination of a  curve in the polar set, 
 with tangents that intersect the given curve (\cite{Lie96} Ch.10, \cite{Lie78}). 
 
 Assume $\eta$ a curve in $\Omega$ and that there is a curve in $\complement \Omega$, such that 
 $d \gamma \cap \eta \neq \emptyset$ (transversal). Note that $d \gamma \bot d \eta$ implies 
 $d \gamma \cap d \eta \neq \emptyset$, that is for every curve $\eta$ In $\Omega$, there is a 
 curve $\gamma \in \complement \Omega$, such that $d \gamma \bot d \eta$.
  
 Assume over $\gamma$ polynomial,  $U$ is harmonic and $U=U_{\times}$, that is a sequential movement. Given
 $\Sigma = \{ U^{\bot} \gamma \quad U \in \mathcal{G} \}$ is planar, that is if for some $U$,
 $U^{\bot}$ gives a planar algebraic geodetic (transversal) in $\Sigma$, that is evolute to $U$,
 we have that $\Sigma$ is algebraic. (cf. strict pseudo convexity (\cite{Oka60})).

 \subsubsection*{Algebraic polar}
 
 Assume $d p^{\bot}=d \gamma$ where $\gamma$ algebraic in the evolute to $p$. 
 Thus given the polar algebraic, we have that $\Sigma_{2}=\{ \gamma \quad X_{U}=0 \}$ 
 is algebraic. Assume $U \gamma$ harmonic and $U$ is algebraic in $(x,y)$, then we have that $p^{\bot} = W \gamma$ for a movement, 
 given $W \gamma + W^{\bot} \gamma=\gamma$ (projective).
 
 Assume $\Sigma=\{ X(f)=0 \}$ and $\Sigma_{0}=\{ \delta X(f)=0 \}$. Then we have that on 
 $\Sigma \cap \Sigma_{0}$, when
 $U^{\bot} f(z)=0$ and $(U^{\bot})^{\diamondsuit}f(z_{1})=0$, that
 $z \sim -1/\overline{z_{1}}$ when $\mid z \mid=1$ (conjugated by involution).
 Assume now $U^{*} \rightarrow \overline{U}^{\diamondsuit}$ projective.

 \newtheorem{polarweak}[adj]{Lemma The polar relative norm}
\begin{polarweak}
 Assume for $f,g \in \mathcal{D}_{L^{1}}$, $<f,g>=(f,\overline{g})$
 and $<f,U^{*} \widehat{g}>=(f,\overline{\widehat{Ug}}) \sim <f,U^{\diamondsuit} 
 \widehat{g}>=0$. Thus, when $f \bot U^{\diamondsuit} \widehat{g}$ the conjugated 
 image can be defined in a closed space.  
 Note that $\parallel U^{\bot} \widehat{g} \parallel_{L^{2}}=\parallel U g \parallel_{L^{2}}=0$
 where $\parallel U^{\bot} \widehat{g} \parallel_{L^{1}}=0$ defines the polar.
 \end{polarweak}

 Assume $\mathcal{G}$ very regular, that is $\exists V \in \mathcal{G}$, such that $(V,V^{\bot})$
 applied on $(f, \widehat{f})$ is regular analytic. Assume $d U / d V=\alpha$, where 
 $\alpha=e^{\phi}$ with $\phi$ pseudo convex in the finite plane. Given $(V,V^{\bot})$ harmonic, 
 we can choose $\phi$ real. Note that $d V \rightarrow 0$
 iff $\frac{1}{\alpha} d U \rightarrow 0$. Given $1/\alpha$ polynomial and real, we have that $dU$ 
 is harmonic outside an algebraic set.

 Given $\log Uf \in L^{1}$, we can assume $U^{\blacktriangle} f \in L^{1}$ (\cite{Dahn18}).
 Assume $X(f)$ defines $U$ according to $(\xi,\eta)$. When $U$ is harmonic, 
 we have $\eta_{x} - \xi_{y} \equiv 0$. Assume $X^{\bot}$ defined by $\overline{U^{\diamondsuit}} \sim U^{\triangledown}$, according to $(\xi^{\bot},\eta^{\bot})$ 
 and that we have $X + X^{\bot}=I$, that is $I(f) \sim \overline{\delta}f$.
 Thus $\xi_{y}=\eta_{x}$ where given $U$ harmonic, $-\eta_{y}=\xi_{x}$,
 that is $\xi(x,y)$ and $\eta(x,y)$ are both symmetric in $x,y$. Conversely given $\xi,\eta$ 
 symmetric in $x,y$ with $\frac{d}{d x}(\eta - \xi)=\frac{d}{d y} (\eta + \xi)=0$, we have 
 $X + X^{\bot}=I$. Assume now $d U^{\bot} / d U = \alpha$ where $\alpha=e^{\phi}$, with $\phi$ 
 pseudo convex outside a small neighbourhood of the boundary.

\subsection{Cylindrical domains}

  A cylindrical domain in this article, is given by $(U_{1},-U_{1}^{\bot}) \simeq (U_{1},-U_{2})$, 
  where we assume $U_{1},U_{2}$ harmonic. When a movement changes character $U \rightarrow I \rightarrow V$,
  we assume a neighbourhood of $I$ is cylindrical. Let $\frac{d U}{d U_{1}}=\alpha$ and $\frac{d U_{1}^{\bot}}{d U_{1}}=\beta$. 
  We assume $\alpha$ locally algebraic, when $U^{\blacktriangle} \rightarrow 0$. Given $\alpha=e^{\phi}$ with 
  $\phi$ pseudo convex, we have $\alpha=1$ iff $\phi=0$ and a neighbourhood where 
  $U \rightarrow I$ is given by $\phi \geq 0$.
  Given $d U$ corresponding to $U$ absolute continuous,has a global base, we have that 
  the transformation that corresponds to $U^{\blacktriangle}$ has a global base. This representation 
  has a maximal domain for absolute continuity in $(x,y)$ (\cite{Malgrange59}). We assume that when the movement changes character, 
  it does not simultaneously change orientation, that is for instance $\alpha \geq 0$.
  Further $d U / d U_{2} \simeq \alpha / \beta$ is a polynomial, given $\beta \mid \alpha$. 
  Consider $U_{1} \beta \simeq \frac{1}{2} \frac{d U_{1}^{\bot 2}}{d U_{1}}$, we assume $U_{1}$
  preserves polynomials. Note that 
  $(\frac{d U}{d U_{1}}, -\frac{d U}{d U_{1}^{\bot}}) \simeq (\alpha,\alpha / \beta)$ and 
  then using transversality, we must assume $\beta \neq const$, 
  why it follows that $\alpha / \beta \prec \prec \alpha$.
  The continuation $U \rightarrow U^{\bot}$ is through 
  $U \rightarrow I \rightarrow U^{\bot}$ $(=U^{\triangle}$) and $(I,0) \rightarrow (U_{1},-U_{2})$

 \newtheorem{lemma}{Lemma the boundary}[section]
 \begin{lemma}
 Given collar point, we have $U_{1}=U_{2}$ on an interval implies $U_{1}=U_{2}=I$. 
 If the movements are analytic and non trivial over $f$, we have that the interval is a point
 \end{lemma}

 Note however that given $U_{1}^{*} \simeq U_{2}$, where $U_{1}$ is given by $(\xi,\eta)$
 constants, we have that $U_{1}^{\diamondsuit}$ is given by $(-\eta,\xi)$ constants. 
 Assume the continuation to $L^{1}$, such that $\mid \xi^{\triangle} / \eta^{\triangle} \mid$ 
 constants, then the completion may include rotation.  
 
 Assume $U^{\triangle}f(\zeta)=f(\zeta_{T})$. Define $\Delta=\{ f(\zeta_{T}) - f(\zeta) =0 \}$. 
 Thus, $U^{\triangle} f - f=0$ over $\Delta$. A necessary condition for hypoelliptic is thus that $(U^{\triangle}-I)f=0$ 
 implies $\zeta_{T}=0$. In particular $\Delta$ defines a ``transversal'' 
 $\subset$ the polar. 
 
 Let the boundary $\Gamma$ be points where the movement changes character, that is 
 $\Gamma(f) = \{ U^{\blacktriangle}f=0 \quad U \in \mathcal{G} \}$.
 Given that we consider the equations in $\mathcal{D}_{L^{1}}'$, we have $<\frac{\delta}{\delta x}(UI),g>=
 <UI, \frac{\delta}{\delta x} g>$. 
 
 \newtheorem{Weyl}[adj]{Lemma Weyl}
 \begin{Weyl}
 When we consider $I_{A}(f)=<X,\frac{\delta f}{\delta x}> + <Y \frac{\delta f}{\delta y}>=0$
 we have $<X_{x} + Y_{y},f>=0$. Let $A^{\diamondsuit}(f)=-Y \frac{\delta f}{\delta x} + X \frac{\delta f}{\delta y}$.
 Then  $I_{A^{\diamondsuit}}(f)=-<Y_{x} - X_{y},f>=0$, given $f \bot (M,W)$ (cf. the projection method \cite{Dahn13})
 \end{Weyl}
 
 Note that according to Weyl, if both the scalar products are $0$ in $L^{2}$, 
 we have $f \in C^{1}$,that is continued to $C^{\infty}$, the orthogonal is in $\mathcal{D}_{L^{2}}$.
 Note that if $T=U^{\blacktriangle}I$ is a distribution in $\mathcal{D}_{L^{1}}'$, we have over the 
 lineality, that $T$ has an infinite zero, that is $<(\delta /\delta x)^{\alpha}T, g>=0$ 
 over $\Delta(g)$. Given that $U^{\blacktriangle}$ can be defined as a measure over $(\delta / \delta x)^{\alpha} g$ 
 $\forall \alpha$, that is given $U^{\blacktriangle}$ is nuclear over all derivatives, we must have that
 $T$ is a $C^{\infty}$ function in the ordinary sense (\cite{Schwartz66}). 
 
 \emph{Given $C$ defined by $U \rightarrow U^{\blacktriangle}$ and denote with $\tilde{C}=S$ a compact, 
 non-discrete connected set, corresponding to $U^{\triangle}$. Then, we have that $U^{\triangle}$ does not necessarily preserve hypoellipticity 
 over $\tilde{C}$.}
 
 The condition $d U^{\triangle}/ d U_{1}=e^{\varphi^{\triangle}}$, 
 where $\varphi^{\triangle}$ is pseudo convex (compact sub level surfaces) implies that $\tilde{C}$ as above
 has finite Dirichlet integral, when $\xi,\eta \in L^{1}(d U_{1})$. 
  
 Assume $U$ represents a sequential movement, that is relative the movement parameters 
 it is given by a cube. Given a proper mapping $\psi$, $U$ can be mapped on to a  ''sliding movement``, 
 where $\psi$ maps collar point on a point in $\infty$. Assume $d U_{j} / d U_{1} = p$, where $p$ is a polynomial 
 locally, we then have, given $U_{j}$ absolute continuous,that $U_{1}=U_{j}$ on an at most algebraic set.
 Note that $d U_{j}^{2} / d U_{1}  \sim U_{j} p$, is not necessarily a polynomial. 
 When polynomials are preserved, given for instance $U_{1}^{2}(t) \sim U_{1}(2 t)$, we can consider $U_{j}$ as algebraic in 
 $U_{1}$. The condition is dependent of division.

 Given the movement is sequential, we represent the movement as with constant dimension on compact 
 sets. The movement parameter interval, is divided into segments of constant dimension 
 for movement.

 When ${}^{t} V^{\blacktriangle}$ is locally 1-1, 
 we see that the movement can be factorized into sequential movements. Note that the 
 collar point condition for ${}^{t} V$ does not imply the same property for $V$, considered in $\mathcal{D}_{L^{1}}'$. 
In the theory of micro local analysis, it is 
sufficient to consider translation domains and otherwise $d U_{j} / d U_{1} = \alpha$ and so on. 
We can complete $\alpha=e^{\phi}$ to $\tilde{\phi}$ harmonic, for instance we assume 
$\log (U \rightarrow U^{\triangle})$ a translation domain.
Polar sets are discussed as defect spaces. 

When $U$ preserves polynomials, that is $\frac{d U}{d U_{1}}=\alpha$ polynomial and
$\frac{d U^{2}}{d U_{1}}=\alpha'$ with $\alpha' \sim \alpha$, then $U \rightarrow U^{2}$
preserves character locally. Note that when $U$ preserves convexity, we have $U - I > 0$
locally. When further $U \rightarrow {}^{t} U$ preserves character, we have $UF {}^{t} (UF) \simeq 
U^{2} F {}^{t} F \simeq U F {}^{t} F$.
When $(U,-U^{\bot}) \rightarrow {}^{t} (U,-U^{\bot})$ preserves character, for instance
$\phi_{1}(U)=-U^{\bot}$ 1-1 with $\phi_{1}^{2}=Id$ and when $U^{*} \rightarrow U^{\bot}$
algebraic, the defect indexes are zero, that is $(U,-U^{\bot})$ is global.

 \newtheorem{maxrank}[adj]{Proposition Global boundary condition}
 \begin{maxrank}
 Given surjectivity for $U$, we have a factorization of ${}^{t} U$ into sequential movements.
 \end{maxrank}

 Assume $W \sim \log S$ with $S \sim (U,-U^{\bot})$ and $V=\frac{S - iI}{S + i I}$. 
 Note that $S + i I$ is interpreted as $(U,-U^{\bot \blacktriangle})$ and $S-iI$ as $ (U,-(-U^{\bot})^{\blacktriangle})$
 We can assume $W \sim 0$ over the diagonal to $\mathcal{G} \times \mathcal{G}^{\bot}$, that is
 we assume that $S \neq (U,U)$. When $U=U_{1}$, we have for the polar condition $U_{1}^{\bot} \sim -U_{1}^{\bot}$, 
 further the lineality is standard complexified. The same condition
 for $U_{2}$ is interpreted such that the rotation is independent of orientation, for instance 
 $\log \mid U_{2} f \mid \in L^{1}$.
 
 Assume $\Gamma=\{ U^{\bot} f=0 \}$.  When 
 $\Gamma \sim \overline{\Gamma}^{\diamondsuit}$ and when the completion to $L^{1}$ is algebraic, 
 we have that $\Gamma \sim \Gamma^{\bot}$.
 Assume $U f = \int f d \mu$, given a global base for $d \mu$, when $U$ is absolute continuous,we have
 a global base for $U^{\blacktriangle}$.  A regular approximation $\gamma_{j} \in \Gamma$, is such that $\frac{d \gamma_{j}}{d t} \neq 0$,
 that is $d y/ d x \neq 0$ and $\eta / \xi \neq 0$, when $\xi$ or $\eta \neq 0$. Consider $U g * f = g * V f$, 
 where $\{ U=0 \} \rightarrow \{ V=0 \}$ is projective, we then have $U \sim V$. When for instance ${}^{t} U \sim U^{\bot}$,
 we obviously have that $\Gamma \sim \Gamma^{\bot}$.
 
 Consider the continuation of $(U,U^{\bot})$, when
 $\{ d U =0 \} \rightarrow \{ d U^{\bot}=0 \}$ is projective, we then have that the defect index
 are equal. A contractible domain corresponds to defect index equal and zero.
 Assume $U F(\gamma) \sim F({}^{t} U \gamma) \rightarrow 0$
 implies ${}^{t} U \gamma \rightarrow 0$. When ${}^{t} U$ injective, we have $\gamma=\gamma_{0}$
 that is $\Gamma \sim 0$ and we have a contractible domain.
 
  \newtheorem{Cousin}[adj]{Proposition  The Cousin problem}
 \begin{Cousin}
 Assume $\mathcal{G}$ has a harmonic generator $U_{1}$ with $U_{1} \rightarrow U^{\bot}_{1}$ is projective,
 such that we have a continuous solution to the second Cousin problem. 
  Then there is a global base for $\mathcal{G}$. For a sufficiently fine 
 division, we assume that $U \in \mathcal{G}$ can be approximated by a sequential movement, 
 arbitrarily close to $U$. Given $U^{\bot}=U^{*}$, $U^{\bot} g=0$ defines  the polar.
 \end{Cousin}
 
 Consider $U$ 
 acting on an ideal $(I)$ and define $(J)$ as the sub ideal where $U$ is real. 
 Thus, if $U=U_{1} + i U_{2}$, we have $(J) \subset \mbox{ ker }U_{2}$.
 In particular, when $U_{2}=U_{1}^{\bot}$ we have that $(J) \subset$ the polar (\cite{Dahn18})
 
 Assume a global base for $\mathcal{G}$, that is $U g=\Sigma U_{j} g$, where $U_{j} \in \mathcal{G}_{j}$
 and $g \in (I)$. 
 We consider monotropy as $\mid U - U^{\times} \mid < \epsilon$ on a non trivial parameter interval.
 Note that a global base  does not imply HE.

 \subsection{Propagation direction}
 
 The condition on collar point assumes the propagation directions the same and is dependent
 on the division intervals to  $\mathcal{G}$. The dimension for $\mathcal{G}$ can vary with $y'$. 
 Assume $\delta(V_{1} f)=U_{1} f$ with $V_{1} f = \int f d \mu_{1}$, then we have given 
 $V_{2} \delta f= \delta (V_{1} f)$, that $<\delta f, d \mu_{2}>=U_{1} f$, and given 
 $\delta \log f=\alpha$, that $<\alpha f,d \mu_{2}>=U_{1} f$.
 When $U_{1} f$ is of negative type, we have $V_{1} f$ is of negative type, 
 that is $V_{1} f$ is regularizing 
 when $U_{1} f$ is regularizing.

 Assume $Y/X$ denotes propagation direction for $Z$ and $-X/Y$  propagation direction for $Z^{\diamondsuit}$.
 Assume further $Y/X \rightarrow -X/Y$ projective, that is preserves exactness.
 When both $Z,Z^{\diamondsuit}$ are harmonic and exact, we have that $Z$ is analytic. 
 The projective mapping implies reflexivity.  Any movement in the plane can,
 according to Lie, be given by $y'' - w(x,y,y')=0$ (\cite{Lie91}). 
 
 Given the movement is harmonic, we have that $\nu \rightarrow \vartheta$ projective.
 Note that the roots to $V,V^{\diamondsuit}$
 are conjugated (through involution). Given $d V,d V^{\diamondsuit}$ are closed, we have that $d V$ is harmonic (\cite{AhlforsSario60}).
 Assume $X(f)$ is defined by $\xi,\eta$ such that $\eta_{x} - \xi_{y} \equiv 0$ (harmonic). 
 In particular when $\xi,\eta$ are constants (translation),  we always have 
 a harmonic movement to start from, for a very regular Lie algebra.
 
 Note that if $U \rightarrow U^{*} \rightarrow U^{\triangle}$ is harmonic,
 we can find a $V \rightarrow \overline{V}^{\diamondsuit} \rightarrow V^{\triangledown}$
 such that, where $X_{U}(f)=X_{V}(f)=0$, we have $V^{\triangledown} \simeq U^{\triangle}$
 
  \subsection{Nuclear movement}
 The polar is dependent on nuclearity.
 Assume $K$ a kernel for an integral operator on $\mathcal{D}_{L^{1}}$, where we have in the weak sense 
 $I_{\widehat{K}}(\phi)=I_{K}(\widehat{\phi})$. Given $\mid <\widehat{K},\psi \otimes \phi> \mid \leq C \parallel \psi \parallel \parallel \phi \parallel$,
 we see that $\parallel \widehat{\phi} \parallel$ defines a polar set, this does not imply $\widehat{K}=0$
 on this set.
 
 A normal surface has finite Dirichlet integral. Assume $\frac{\delta {}^{t} U}{\delta y}=\overline{\eta}$ (or $\overline{Y}$),
 such that $\int d U = \int \frac{d U}{d t} d t = \int \{ F,{}^{t} U \} d t$ finite, where we 
 assume the movement dependent on one parameter. We thus assume $Y \overline{\eta} + X \overline{\xi} \sim \mid X \mid^{2} + \mid Y \mid^{2}$.
 Assume further $F^{\triangledown}F(x,y)=(F(\tilde{x},y),M(\tilde{\nu_{1}},\vartheta_{1}))$ (\cite{Dahn18}). 
 Through the condition
 $F^{\triangledown} F \sim F F^{\triangledown}$, we see that $F(\tilde{x},y) \sim F(x,\tilde{y})$. '
 
 Assume $C$ the set where $U \rightarrow U^{\bot}$ and $\tilde{C}$ the set corresponding to the movement completed  
 to $L^{1}$, that is $U=U^{\bot}$ ($U^{\triangle},U^{\triangledown},U^{\blacktriangle}$) over a compact set $\tilde{C}$. Assume in a cylindrical neighbourhood 
 of $\tilde{C}$, there is a
 sequential movement $U$, dependent of division.  Assume $\mid d U \mid=d U(\nu,\vartheta)$, 
 then on $\tilde{C}$ we have that $\nu=\vartheta$.

 \subsection{The extended plane}
 
 Consider $U (x y(x))$, then we have that $\frac{\delta U (x,xy)}{\delta x}=\xi y$
 and $\frac{\delta U (x,xy)}{\delta y}=\eta x$, that is if $U$ is translation, the same movement on
 $(x,xy)$ is rotation and given $X \rightarrow X^{\bot}$ defines a cylindrical domain, 
 we have that $(x,y) \rightarrow (x,xy)$ corresponds to $U \rightarrow U^{\bot}$. 
 Assume for $U,V \in \mathcal{G}$, we have $dU / d V=\alpha$.
 
 \newtheorem{ahnlich}[adj]{Lemma The \"ahnlich transform is proper}
 \begin{ahnlich}
 Assume $C$ corresponds to constant surfaces for $\alpha$. Given $\psi : $ discrete 
 $\rightarrow $ compact (cf almost 1-1 coverings) a change of character is of the form 
 $\psi^{*} \alpha=1$. Note in particular that $\{ \alpha < \lambda \}$ compact is implied by 
 $\{ \alpha(x,\frac{y}{x}) < \lambda \}$ compact, assuming $y/x < y$ in $\infty$. That is 
 $(x,y) \rightarrow (x,y/x)$ acts as a proper mapping.
 \end{ahnlich}
 
 Assume $V$
 a continuous movement acting on $(x,\frac{y}{x})$ close to $U$ analytic acting on $(x,y)$, then
 a discrete boundary to $U$ corresponds to a compact boundary to $V$.
 When we only have compact sub level surfaces for the extended system, the system does not preserve
 constant value in the $\infty$ and $0$ simultaneously in both variables. Given $M$, the symbol in $\nu,\vartheta$, 
 symmetric and not algebraic (semialgebraic), there is possibility of presence of spiral.
 The condition on collar point is only assumed for $x,y$, not the extended plane.
 
 Assume $I_{M}(\phi)=\int M(\nu,\vartheta) \phi(\vartheta) d \vartheta$ and that 
 $U^{\triangledown}_{1} M=0$ implies $\mid \nu \mid \leq 1$, that is $\mid \vartheta \mid \geq 1$. 
 Then we have that $U=U_{1}$ implies compact
 sub level surfaces in $\nu$, that is $U^{\triangledown} M \leq \lambda$ implies 
 $\mid \nu \mid \leq 1$. Simultaneously unbounded sub level surfaces in $\vartheta$. 
 
 \newtheorem{regularM}[adj]{Proposition Regularity in extended domain}
 \begin{regularM}
 Thus, if ${}^{t} I_{M}$ regular in $\nu$ and $\vartheta$,
 we have that $M$, if non-trivial, is not symmetric. That is, if $(F,M)$ is regular in the extended plane, 
 then $M$ is not symmetric.
 \end{regularM}

 Note however, that $M$ can have compact sub level surfaces in $\zeta$, when $M$ is symmetric in 
 $\nu,\vartheta$. Consider $F(x,y,\frac{d y}{d x}, \frac{d x}{d y})$. Given 
 $\phi (\frac{y}{x})=\frac{x}{y}$ projective, we can shorter write $F(x,y,\frac{d y}{ d x})$. 
 Assume $K$ a very regular Schwartz kernel, in this case $(K,M)$ means that $M$ is singular 
 on $x=y \sim (z,z)$. Consider $K(d x,\frac{d U y}{d t}) \sim K(d x,d U)$, 
 when $x,y$ fix. We write $\int d K \sim \int K(dx, d U) K(d U^{\bot}, d y)$.

 Assume now $\gamma$ a curve, such that there is a $\tilde{\gamma} \in \mathcal{K}$,  
 a class of curves, such that $\mid \tilde{\gamma} - \gamma \mid \leq \epsilon$. 
 Assume $\mathcal{K}$ a boundary class (involution, reflection) 
 and consider $\Omega_{\gamma}$  such that $\mid \gamma \mid \leq \mid  \tilde{\gamma} \mid$ (length of curve) 
 where $\tilde{\gamma} \in \mathcal{K}$, then $\Omega_{\gamma}$ includes geodetics and is for instance a spiral 
 region. Assume now $\mathcal{K}$ a class of curves
 for which the continuation preserve algebraicity. Given existence of$\tilde{\gamma}_{1},\tilde{\gamma}_{2} \in \mathcal{K}$
 such that $\mid \tilde{\gamma}_{1} \mid \leq \mid \gamma \mid \leq \mid \tilde{\gamma}_{2} \mid$, we have that for a spiral 
 $\gamma$, $\gamma \notin \mathcal{K}$. Consider now Iversen's model (\cite{Dahn13}), that is 
 $0 < \mid \log \tilde{\gamma} - \log \gamma \mid < \epsilon$,
 and assume that the inequality can be continued algebraically over 0. Thus, 
 $\log \tilde{\gamma} - \epsilon < \log \gamma < \log \tilde{\gamma} + \epsilon$.
 Given the condition for algebraic continuation, $\gamma \in \mathcal{K}$.
 The condition is sufficient for HE.

 \newtheorem{symplorto}[adj]{Lemma Symplectic orthogonal for the graph}
 \begin{symplorto}
 Assume $\mid d \mu \mid = d \mu(\nu,\vartheta)$, where $d \mu^{\triangle}=p d \mu^{*}$,
 for a polynomial $p$, is the completion to symplecticity. For the extended system, 
 we have that $<d v, d \mu^{\triangle}>=0$
 $\forall d v \in T_{p} \mathcal{G}$ implies $d \mu^{*}=0$. Given symplecticity, 
 assume for $\forall U \in \mathcal{G}$ that $<d w, d (U,U^{\bot})>=0$. Assume 
 $d w \sim d (V,V^{\bot})$, where $V \in \mathcal{G}$ and absolute continuous. Given $d (V,V^{\bot})$ integrable, 
 we have that $(V,V^{\bot})^{\blacktriangle}=\{ 0 \}$
  \end{symplorto}

 \subsection{Global concepts}

  Concerning involution: Assume $\{ {}^{t} V, U \}=0$, that is $\frac{\delta {}^{t} V}{\delta x} \frac{\delta  U}{\delta y}= 
 \frac{\delta {}^{t} V}{\delta y} \frac{\delta U}{\delta x}$. Then we have that ${}^{t} V$ and $U$ has 
 the same propagation direction. Given $\frac{\delta {}^{t} V}{\delta x}=\frac{\delta V}{\delta y}$ 
 and $\frac{\delta {}^{t} V}{\delta y}=-\frac{\delta V}{\delta x}$,
 we have that $(U,V)$ (that is $U + i V$) can be considered as a analytic movement. Note that 
 given $d U = \alpha d V^{\bot}$, we have that $\alpha=const$ over an involutive set (cf. evolute).
 Thus given $U \in \mathcal{G}$, we assume that we have existence of $V \in \mathcal{G}^{\bot}$ 
 such that $U \simeq V^{\bot}$ over an involutive set.
 Consider now the completion $V^{\triangle}$, that is we assume $(U,V^{\triangle})$ closed. 
 The notation is obviously improper for functional, but can be motivated for $\mathcal{D}_{L^{1}}'$.
 
 Note that when $d w$ is exact and $d w^{\bot} \gamma \rightarrow 0$ implies $\gamma \rightarrow 0$, 
 then we have that for the corresponding domains $\Omega_{w} \cap \Omega_{w^{\bot}} \sim \{ 0 \}$. 
 Sufficient for this is that $dw , d w^{\bot}$ are exact (cf. harmonic conjugation).
 
 \newtheorem{1_polar}[adj]{Definition 1-polar for the graph (L. Schwartz)} 
 \begin{1_polar}
 If relative $L^{1}$- norm, $U^{*} \widehat{f}=0$ and $\widehat{(I-U) f}=0$, we have that $\widehat{f}=0$. Further, when 
 $(U f,U^{*} \widehat{f}) \rightarrow 0$ and when $U \rightarrow I$ regularly, 
 $f \rightarrow 0$. We define the polar, so that $U^{*} \widehat{f}=0$
 for $\widehat{f} \neq 0$. 
 \end{1_polar}
 
 Concerning Oka's property (\cite{Oka60}, \cite{Range}): Assume existence of $d U_{j}$ analytic over $\Omega_{j}$ and 
 $d U_{ij}=\frac{d U_{j}}{d U_{i}}=\alpha_{ij}$ analytic and not zero on $\Omega_{i} \cap \Omega_{j}$. 
 A global base means that we have existence of $d U$ such that $\frac{d U}{d U_{j}}$ is continuous and 
 not zero on  every $\Omega_{j}$.
 Oka's property for $e^{\phi_{i}}$, means $\phi_{ij}=-\phi_{ji}$ on $\Omega_{i} \cap \Omega_{j}$.
 further $\phi_{ij} + \phi_{jk} + \phi_{ki}=0$ on $\Omega_{i} \cap \Omega_{j} \cap \Omega_{k}$.
 Thus we have existence of $\phi_{j}$ continuous on $\Omega_{j}$ (cf. pseudo convexity) such that 
 $e^{\phi_{ij}}=e^{\phi_{i} - \phi_{j}}$. Thus $\alpha_{j}$ can be chosen to give a global 
 pseudo base for the movements, as above.  
 
 Given Oka's property for analytic movements, let $(I_{2})=$ $\{ f \in H(\Omega_{2}) \quad X_{2}(f)=0 \}$
 and on $\Omega_{i} \cap \Omega_{j}$, $X_{1}(f)=X_{2}(f)=0$. When we refer to a sequential
 movement, a movement is divided using $C$ into segments where the movement is determined, 
 $\Sigma_{i}=\{ (x,y) \quad U = U_{i} \} = \cup_{ij} \Sigma_{ij}'$.
 In this case, we assume when $\Sigma_{i} \cap \Sigma_{j} \neq 0$, that we have $\frac{d U_{i}}{d U_{j}}=\alpha_{ij} \in H((nbhd C) \backslash C)$.

 \newtheorem{ac}[adj]{Lemma Absolute continuous transformations}
 \begin{ac}
 Starting with a global base for $d U$, we can determine a maximal domain
 for absolute continuity.On this domain, we have that $d U=0$ implies $U^{\blacktriangle}=0$.
 \end{ac}
 
 (\cite{Malgrange59}) Obviously, $U-cI=0$, where we assume $c=1$ (modulo scaling).
 Given $\tilde{C}$ corresponding to $U^{\triangle}$
 and $\psi : C \rightarrow \tilde{C}$ continuous and proper, if the continuation is analytic, 
 we have that $\tilde{C}$ is removable iff we have a global base for the corresponding ideal. 
 In particular, when the continuation is algebraic, we have that $\tilde{C}$ is removable.
 
 \newtheorem{monotropy}[adj]{Lemma Monotropy for $\mathcal{G}$}
 \begin{monotropy}
 Assume $\mathcal{G}$ very regular, in the sense that for every
 symbol $f$, there is a movement $V$, analytic over $f$. Assume ${}^{t} (I-U) \gamma \rightarrow 0$ 
 implies $\gamma \rightarrow \gamma_{0}$, then there is a $V$ such that ${}^{t}(I-V) \gamma \rightarrow 0$ 
 analytic, that is close to $\gamma_{0}$, we have that $\mid {}^{t} (U-V) \gamma \mid \leq \epsilon$, 
 for $\epsilon$ small and positive.
 \end{monotropy}

 \subsection{The wave front set}
 
 G{\aa}rding (\cite{Garding87}) defines the hyperbolic cone $\mathcal{C}(f,a)$ for $f$ in $a$, as the component in the complement
 to the real hyper surface $Pf f(x)=0$, that contains a. The dual cone $K(f,a)$ to $\mathcal{C}(f,a)$ such that
 $x . \mathcal{C} \geq 0$ is a closed and convex propagation cone. The wave front set is a closed semi algebraic
 subset of the propagation cone of co dimension 1.  
 
 Assume $\Sigma=\{(x,y) \quad X(f)=0 \}$, where $X(f)=0$ is corresponding to $U$ on $\Sigma$. Then we have that
 $X(f)(tx,ty)=0$ for $(x,y) \in \Sigma$, that is $U$ can be defined as independent of scaling. In this manner,
 $(\xi,\eta)$ can be regarded as in a dual relationship to $(x,y)$ (with respect to scaling).  Define 
 $\mathcal{C}(f,a)=\complement \Sigma$. Thus we have that $(x,y) \in \mathcal{C}$ implies $t (x,y) \in \mathcal{C}$.
 Given a normed space, it is to define
 the propagation cone, sufficient to consider $(x,y).\mathcal{C}=0$, that is we can consider $\mathcal{C}^{\bot}$
 
 A co dimension one variety is defined by $S(p,x)=0$ and $s_{x} \neq 0$
where $p = s_{x}$ is characteristic. Assume $K(p)$ defined by $g(x,s_{x})=d g(x,s_{x})=0$, $p.dx=0$ and $d S=p dx$,
then we have that $S(p)$ is involutive if $H_{g}=0$ (Hamiltonian).
  We define $K=\{ (\xi,x) \quad x \in K(\xi) \mbox{  iff } K(x) \cap S(\xi) \neq \emptyset \}$.
 If $x(t,z)$ is the uniformization of $u^{\mu}$ up to a certain order in a non-characteristic point,
 we have that $\mbox{ sing supp } u^{\mu}(x) \subset K$ (bicharacteristics such that $K \cap S \neq \emptyset$)
(\cite{Garding64})
 \newtheorem{bikar}[adj]{Definition Bicharacteristics}
 \begin{bikar}
 Assume the symbol $g(z,p)=0$, corresponds to $X(f)=0$ to $U$ and assume $S(\xi,\eta)$ is defined
 by the condition 
 $U^{\bot}f=0$, that is polar points and that characteristics (for $U^{\bot}$) are given by $X(f)=0$. Bicharacteristics corresponds to
 $(\xi,\eta) \sim p$ $\bot (-Y,X)$. Further $\frac{d p}{d z} \sim -\frac{g_{z}}{g_{p}}$. Given $f$ a Hamiltonian
 we have that consequently $g_{p} \sim (-Y,X)$. 
 \end{bikar}

 The condition that we have existence of
 $U \in \mathcal{G}$, where $U_{t} \rightarrow U$, such that $U_{t} f$ is holomorphic in a parameter neighbourhood,
 can be compared with the condition that $\mathcal{G}$ is very regular. In this manner the analyticity is preserved under 
 the uniformization and in particular the algebraicity is locally preserved. The condition $K \cap S \neq \emptyset$
 means that we have existence of a regular approximation of a singularity. Since we assume
 $f \in \mathcal{D}_{L^{1}}$, we can assume the complement to $S$ a translation domain, that is
 1-dimensional. Finally, note that when we consider
 $U F(\gamma)(\zeta)=F({}^{t} U \gamma)(\zeta)=F(\gamma)(\zeta_{T})$ and when $U \rightarrow \zeta_{T}$
 is continuous, we can transport the argument to $\zeta_{T}$.
 
 Consider $(Uf,U^{*} f)$ as a graph in $\mathcal{D} \subset \mathcal{D}_{L^{1}} \subset \mathcal{D}_{L^{2}}$.
 Given $U^{*} f$ can be continued to $L^{1}$ and $L^{2}$, consider $(Uf,U^{\triangle}f)$ and 
 $(Uf,U^{\bot} f)$ respectively.
 Given $U^{\triangle}$  $\rightarrow Uf$ locally 1-1, we have a topological ``monodromy''. 
 For the continuation to $L^{1}$ we consider Cauchy sequences,  $\parallel \phi_{n}^{*} - \phi^{\triangle} \parallel \rightarrow 0$
 why in $L^{1} \cap L^{2}$, we have that $\phi^{\bot} \sim \{ \psi \quad \psi \bot \phi^{\triangle} \}$. 
 For hypoelliptic $f$ we have that $\phi^{*} \sim \phi^{\triangle}$
 and for harmonic movements,  that $\phi^{\bot} \sim \{ \psi \quad \psi \bot \phi^{\triangledown} \}$

\section{Symplecticity}
 
 Assume $F \gamma = e^{\phi^{*} + \phi_{1}}$, where $\phi_{1}$ polynomial. Assume 
 $\Omega_{0}=\{ \phi_{1}=0 \}$ and that we have existence of a domain of holomorphy  
 $\Omega \supset \Omega_{0}$, where $\phi_{1} > 0$ on $\Omega$.
 Assume $P = N(\phi_{1})$ algebraic with 
 $P \subset \Omega$, then we have that over $P$, 
 $e^{\phi^{\triangle}}=e^{\phi^{*}}$ on $P$.
 Given $\phi_{1}$ polynomial, $\{ d^{2} \mid F \gamma \mid \leq C \} \sim \{ \mid F \gamma \mid \leq 1/Q \}$, 
 for a polynomial $Q$, a domain of holomorphy. (cf preserves constant value according to \cite{Cousin95}).

 Assume $A f=0$ and $A^{\diamondsuit}  g=0$, then we have that
 $\big[ X_{f},X^{\diamondsuit} g \big]=\big[ X_{f},X_{h} \big]$ (Lie-bracket), that is we assume $X^{\diamondsuit} g \sim X h$. Thus,
 given $\xi / \eta \sim X /Y$, we have that $\frac{\delta g / \delta x}{\delta g / \delta y}=\frac{\delta h / \delta y}{\delta h / \delta x}$.
 Consider $X(A g)=\lambda A g$. Given $g$ a Hamiltonian,
 we have that $A g=0$, that is $d g \bot d h$ or $\{ h,g \}=0$. Assume $\frac{\delta h}{\delta x}=\frac{\delta U g}{\delta x}$, when 
 we have that $\{ g,Ug \}=0$ that is $d g \bot d U g$ and given $U$ ac, when $\bot$ is symplectic,
 $U^{\blacktriangle} g=0$ (invariant points). 
 Further $\frac{\delta g / \delta x}{\delta g / \delta y} = \frac{\delta h / \delta x}{\delta h / \delta y}$
 that is given $h$ has no symmetry set, we have that the same holds for $f$, over an involutive set. If $h$ has no non-closed 
 extension over an involutive set, we have that the same holds for $f$.

\newtheorem{trans}[adj]{Definition Dual transformations}
\begin{trans}
Consider $\Sigma=\{ X(f)=0 \quad X^{\bot}(\widehat{f})=0 \}$ and the condition 
$X (X^{\bot} (\widehat{f}))=\lambda X^{\bot}(\widehat{f})$.
We define $X^{\bot}(\widehat{f})=X^{*} \frac{\delta \widehat{f}}{\delta x^{*}} + Y^{*} \frac{\delta \widehat{f}}{\delta y^{*}}$,
that is $(\xi^{*},\eta^{*}) \bot (-Y^{*},X^{*})$.
\end{trans}

Given $\Sigma=\{ X(f)=0 \quad X^{\triangle}(\widehat{f})=0 \}$, we can define $\Sigma^{\bot}$ as 
symplectic orthogonal.
Consider the mapping $X_{x} + Y_{y} \rightarrow \widehat{X}_{x^{*}} + (\widehat{Y})_{y^{*}}$.
Assume $Y/X \sim v$ and $v^{*} \sim \widehat{Y} / \widehat{X}$, that is $\widehat{Y} \sim \widehat{V * X}$,
where $\widehat{V}=v^{*}$. Given $X^{\bot}(\widehat{f})=0$ we have that on a contractible domain,
$- v^{*} = \frac{d y^{*}}{d x^{*}}$.  When $v$ is the propagation direction,
then $v^{*}$ is the propagation direction for the system to $\widehat{f}$. 
From the theory on multipliers, $v^{*}$ is a multiplier if $\parallel \mathcal{F}^{-1} (v^{*} \widehat{f}) \parallel \leq \parallel f \parallel$ (\cite{Schechter70}).
The corresponding convolution operator $V$ is invertible if $v^{*}$ is downward bounded.
Consider for this reason $\{ \frac{d y^{*}}{d x^{*}} < \lambda \}$, that are bounded if $(x^{*},y^{*})$ on one side
of a plane, except for a compact set. For instance when $U^{\bot}$ (or $U$) reflection, or
when $U^{\bot}$ convex.
 
\subsubsection{Symplectic completion}
Consider $\Sigma(g)=\{f \quad X(f)=0 \quad \delta X(f)=g \}$, then we can define $g^{\bot}$ as 
an annihilator, that is when $\Sigma$ closed we have that $\{ f \quad <d f,d g^{\bot}>=0 \}=\Sigma(g)$.
Note that $\mathcal{L} \subset \Sigma^{\bot}(g)$ can be defined as bicharacteristics leaves where the movement is defined as nuclear,
that is the movement can be factorized.
Note that when $U^{\blacktriangle}=I-U$, 
we have that $U^{\blacktriangle}$ has isolated zero's, where $U$ has algebraic zero's, corresponding to a closed extension.

Assume $\mathcal{D}_{L^{1}}=(\dot{B}) \bigoplus B_{0}$ and that for $U$ harmonic, $\mathcal{D}_{L^{1}} \ni \phi \rightarrow U \phi \in \mathcal{D}_{L^{1}}$.
Assume $g \bot \dot{B}$, $g \in \mathcal{D}_{L^{1}}$ with support in a neighbourhood of $\infty$. When $<U^{*} \widehat{\phi},g>=
<\widehat{\phi}, {}^{t} V g>$, for all $\phi \in \mathcal{D}_{L^{1}}$, we have $\int g d {}^{t} (U^{*} - V)=0$
and since $V$ can be chosen as harmonic, why ${}^{t} V \rightarrow V$ preserves character,
we can conclude, on sets where the measures are finite, that $U^{*}=V$.
\newtheorem{kompl}[adj]{Lemma Symplectic condition}
\begin{kompl}
Symplecticity means $X_{j}^{\triangle}(f) \bot \widehat{\psi}$, $\forall \widehat{\psi}$ implies
$X_{j}^{\triangle}(f)=0$. Thus $U_{j} \rightarrow U_{j}^{*}$ can be completed to analyticity.
In particular given $U$ harmonic, we have given symplecticity, that $U \rightarrow U^{*}$ defines the movement $U^{\bot}$ uniquely.
\end{kompl}
  Assume $U f = \int f d \mu$ on $\Sigma$ and $d U^{\bot}=(e^{\phi} -1) d U_{1}$, where 
 $d U_{1}$ harmonic. The polar to $d U^{\bot}$ can when be written $\phi=0$ over a compact. 
 Consider the completion $\dot{B} \rightarrow \mathcal{D}_{L^{1}}$. Assume
 $<U^{\triangle} \widehat{f},g>=<\widehat{f},g>$. Given $U^{\triangle}$ 
 absolute continuous,this is not dependent of division. Given $d U^{\triangle} \widehat{f} = (e^{\phi} -1) d U^{*} \widehat{f}$
 and $\phi \equiv 0$ outside a compact, gives a continuation $\dot{B} \rightarrow \mathcal{D}_{L^{1}}$.
 Given $\delta X_{j}^{\triangle}=\delta (p X_{j}^{\bot})$, we have given $X_{j}^{\bot}$ harmonic and
analytic, that $X_{j}^{\triangle}$ are harmonic and analytic. Otherwise the condition for $\delta X_{j}^{\triangle}(f)=0$
is that $\delta \log X_{j}^{\bot}=- \delta \log p$.
 
 \emph{Consider $X(f) \rightarrow f$, assume $Y d x - X d y=0$, where $X,Y$ are analytic, such that 
 $X=Y=0$ are points. Then, for instance, given $X \neq 0$, we have that $\eta/\xi \sim Y/X$. 
 Given that $X(f)=0$ implies $f=const$, then $f$ is not 
 holomorphic. In the planar case, when $f \in (I)$, the ideal of integral curves, 
 we have that $g \in (I)$ implies $g=const$.} (\cite{Lie91}, Ch. 5, Theorem 3)

 When we introduce monotropy, if $\gamma$ is a cycle corresponding to $ dw= d V^{\bot}$, 
 where $V$ analytic and $\tilde{\gamma}$ is the monotropic consequent to $d V^{\bot}$, 
 then we have that $\tilde{\gamma}$ is not necessarily closed..
 Assume $\int_{\tilde{\gamma}} d w \sim \int_{\gamma} d \tilde{w}$ (\cite{Dahn13}), that is by
 completion of $d w$, we can relate $\tilde{\gamma}$ to a closed contour.
 
 \newtheorem{lacunas}[adj]{Definition}
 \begin{lacunas}
 When $\tilde{C}=\psi(C)$ and
 $d U^{\bot}=0$ on $\gamma \sim 0$. Assume we have existence of $W$ with $W^{\bot}=0$ on $\tilde{C}$ and $d U^{\bot}(0)=0$.
 Then we have that $( W U )^{\bot}=0$ over $\tilde{C}$ and $( U W )^{\bot}=0$ when $U^{\bot} \gamma \in \tilde{C}$.
 Note that given $U_{2} \rightarrow I \rightarrow U_{3}$ in a point, means that $\gamma \sim 0$,
 that is $C$ is discrete when the movement sequential (and analytic), corresponding to $\gamma \sim 0$.
 \end{lacunas}
 
 \subsubsection{A normal model}
 
 The lifting principle can be solved over an algebraic polyhedron (\cite{Oka60}). For an analytic polyhedron, we consider
 a normal model $\Sigma$ (ramified domain), where the lifting principle can be solved. 

 Assume $f \in \mathcal{D}_{L^{1}}$ and $g$ analytic, with $\{ f,g \}=0$. Given monotropy, where
 $(U^{\diamondsuit})^{\bot}$ is considered as an analytic function over $d f$, then we have existence, given
 $(U^{\diamondsuit})^{\bot} d f=0$, of $d g$ analytic such that $(U^{\diamondsuit})^{\bot} d g=0$
 and $\mid d f- d g \mid \leq \epsilon$. Alternatively, let $z$ be parameter for movement corresponding to $(U^{\diamondsuit})^{\bot}$.
 Given $g$ analytic in $z$ and $f$ analytic in $z_{1}$, parameter to $(U^{\diamondsuit}_{1})^{\bot}$,
 then monotropy means that $\mid z- z_{1} \mid \leq \epsilon$, that is the movement relative $U$ can be approximated
 arbitrarily close with the movement relative $U_{1}$..
 
 Assume $U$ is approximated by $U_{\times}$, a sequential movement and $\tilde{\psi} U_{\times} = U_{\times} \psi$.
 Thus, $U^{\bot}=0$ on $\tilde{C}$ and $U_{\times}^{\bot}=0$ on $C$. 
 Assume $(\Gamma)$ is given by $\mid \gamma \mid < 1$
 and $d U^{\bot}$ analytic and finite on $(\Gamma)$. Given $d U^{\bot}=0$ on a segment of $\Gamma$
 we have that $d U^{\bot} \equiv 0$ on $(\Gamma)$ (\cite{Collingwood66}, Theorem of F. and M. Riesz), that is on en set of positive measure. For instance
 $\tilde{\psi} U_{\times}^{\bot}=0$. We assume here $\Gamma$ closed in the plane, that is a pluricomplex definition.
 
 \newtheorem{Dirichlet}[adj]{Propostition The Dirichlet integral}
 \begin{Dirichlet}
 
 A normal surface has finite Dirichlet integral. Assume $\Sigma = \{ \gamma \quad D_{\Sigma}(\gamma) < \infty \}$.
 Given the domain contractible, we have where $U$ is analytic, $\int d U(\gamma)=0$ when $\Sigma \ni \gamma \rightarrow 0$
 regularly.  Consider $\frac{d U}{d t}=\xi \frac{d x}{d t} + \eta \frac{d y}{d t}$
 given $\mid \frac{d U^{\diamondsuit}}{d t} \mid=0$ and $0 < \mid Y/X \mid$ finite, we have $\mid \frac{d U}{d t} \mid \sim \mid X \mid^{2} + \mid Y \mid^{2}$, that is $\int_{\Sigma} \mid d U \mid(f)$
 finite iff $\Sigma$ has finite Dirichlet integral. Further, if $d U$ is harmonic, we have a finite D-integral
 iff $\int d U (f)$ is finite.  
 \end{Dirichlet}
 
 Note that $\overline{d U^{\diamondsuit}}$ has coefficients $\overline{(-\eta,\xi)}$
 which in this context is $-(\eta,\xi)$.
 Assume $(X,Y) \in \mathcal{D}_{L^{1}}$ we can then assume $(\widehat{X},\widehat{Y}) \rightarrow 0$
 when $x,y \rightarrow \infty$. We assume absence of essential singularities in the $\infty$ (\cite{Cousin95}).
 For instance, we can assume $\widehat{X},\widehat{Y}$ preserves constant value in $\infty$. 
 Assume $X=\widehat{X_{1}},Y=\widehat{Y_{1}}$, that is $\widehat{Y_{1}}=\widehat{\rho} \widehat{X_{1}}$.
 when we have that $Y/X \sim \widehat{\rho}$. Given $\widehat{\rho}$ denotes the propagation direction,
 we have that $\widehat{W}= \widehat{Y_{x} - X_{y}} \sim (-y + \widehat{\rho} x) \widehat{X}$ and
 $\widehat{M} \sim (x + \widehat{\rho} y) \widehat{X}$. Further, we have that $\frac{x + \widehat{\rho} y}{-y + \widehat{\rho} x}$
 changes sign in 0 in the same manner as $-\frac{x}{y} - \widehat{\rho}$. Note that when $\widehat{\rho}$ is constant,
 we have that the quotient is $\sim$ the propagation direction for the spiral. In this case we have thus $\nu,\vartheta$ 
 are constants, that is the spiral behaves as the degenerate case for dynamical systems. 
 
 Given that the polar can be defined by $\mathcal{G}$, we can assume finite D-integral.

 \section{Polar sets}
 L. Schwartz (\cite{Schwartz54}) discusses trace functions, $\mathcal{E}_{L^{2}}^{m}(\complement \Omega) \simeq \mathcal{D}_{L^{2}}^{m}(\complement \Omega)$
 iff $\gamma^{r} f=0$ on $\Omega$ for $m - \mid r \mid \geq (\frac{n-p}{2} + 1)$. When 
 $D^{s}$  denotes the tangential derivative, we have that $D^{s} \gamma^{r} f = \gamma^{r} (D^{s} f)$. 
 Assume $V_{\lambda}=\{ f < \lambda \} = \complement \Omega$. Let $f_{\lambda}= f \mid_{V_{\lambda}}$
 then we have that a sufficient condition for $\gamma(f_{\lambda})=0$ on $\Omega$ is that $f$ HE. 

 Assume $\mbox{ Im }f / \mbox{ Re }f \rightarrow 0=e^{\phi} \rightarrow 0$ in $\infty$, 
 that is $e^{\phi} \in \dot{B}$. When $\mathcal{D}_{L^{1}} \ni f \sim \mbox{ Re }f$ and $X(f)=0$, we have existence of 
 $g \in H \cap \mathcal{D}_{L^{1}}$ such that $\Delta g=0$.
 Thus if we assume $\delta X(g)=X(g)=0$, $g$ can be used to continue the movement to 
 $\mathcal{D}_{L^{2}}$. 
 
 Assume $\Sigma=\{ X(f)=0 \}$, $Uf$ not identically constant and $\Sigma_{0}=\{ \delta X(f)=0 \}$, 
 on $\Sigma_{0}$ we have that
 $M(Uf)=Uf$ (arithmetic mean) why $\Sigma_{0} \cap \Sigma$ are isolated points, for harmonic movements.
 Assume $f \in \mathcal{D}_{L^{1}}$ and $\varphi \in \mathcal{D}$ (or $\mathcal{S}$), 
 such that $f * \varphi \in \mathcal{D}_{L^{2}}$ and $V (f * \varphi)=f * {}^{t} V \varphi$. 
 Given $f \in B_{pp}$ (\cite{Schwartz66}) we have that $M(f * {}^{t} V \varphi)=const$ implies $M(V f)=const$ 
 when $\varphi \rightarrow \delta_{0}$. Note $f$ hypoelliptic implies $1/f \in B_{pp}$. 

 \newtheorem{regmovement}[adj]{Lemma Regularization of movement}
 \begin{regmovement}
 Assume $f \in \mathcal{D}_{L^{1}}$ and $g \in \mathcal{D}_{L^{2}}$ and $g \rightarrow \delta_{0}$. 
 Further, $U \in \mathcal({D}_{L^{1}})'$, then given $U$ ``algebraic'' we can approximate $Uf$ 
 by $Uf * g \sim f * U g \in \mathcal{D}_{L^{2}}$. 
 Thus, $U^{\triangle}$ can be determined in $\mathcal{D}_{L^{2}}'$. 
 \end{regmovement}

 Assume $d U / d U_{1}=\alpha$, where $\alpha$ is a polynomial locally. Then we have that 
 $d U^{j} / d U_{1} \sim U^{j-1} \alpha$, that is if $U$ preserves polynomials (sequential movement),
 collar point is preserved for $U^{j} f$ and the same holds for $1/U^{\blacktriangle}$. 
 Note that if $\widehat{U^{\bot} f}=U^{\bot \bot} \widehat{f}$, 
 given $U^{\bot \bot} \rightarrow U$ 
 continuous, we have that the polar, that is $\widehat{U^{\bot} f}=0$ is preserved, that is 
 the polar is generated by $\mathcal{G}$ in this case. 
 
  \newtheorem{ortho}[adj]{Lemma Orthogonal movement}
 \begin{ortho}
 Assume $\phi \in (I)$ implies $U \phi \in (I)$ and $\phi = U \phi + (U \phi)^{\bot}$
 in $L^{1}$. Then we can obviously define a movement $U^{\bot} \phi^{\bot}=(U \phi)^{\bot}$.
 \end{ortho}

 Assume $U \phi(x)=\phi( g x)$, then we can define $(g x) \rightarrow g^{\bot} x^{*}$ 
 (Legendre) that is orthogonal with respect to the circle. Given $U \phi$ analytic, we have 
 that the mapping $U \rightarrow g$ is continuous. When the movement is considered
 in a weak sense, we can consider $H' \rightarrow Exp$. In this case it is necessary to 
 put the condition, ${}^{t} U$ preserves polynomials, since $Exp_{L^{1}}$ does not have an 
 algebraic base.
 
 \newtheorem{lineality}[adj]{Lemma The lineality as a polar set}
 \begin{lineality}
 The lineality to a symbol in $\mathcal{D}_{L^{1}}$, is a polar set.
  \end{lineality}

 Assume $-U_{1}^{\bot} f=0$ (translation and $U^{\bot}_{1}=U^{\blacktriangle}_{1}$) on a set of positive measure (a line), given isolated singularities
 for instance $f \in \mathcal{D}_{L^{1}}$, then we have that $U \rightarrow \zeta_{T}$ is mapped 
 on translation. Thus, we have that $f(\zeta_{T}) \equiv f(\zeta)$ on a line $L$. Conversely 
 given $L \subset \Delta$ (lineality) we have that $U_{1}^{\bot} f \equiv 0$ on a curve in 
 $(x,y)$, which implies that $f$ is not HE. $U_{1}$ considered
 as a distribution, is thus orthogonal (relative $f$) on L, that is $L$ can be seen as a polar.

 \newtheorem{parametrix}[adj]{Lemma The parametrix kernel represents a polar}
 \begin{parametrix}
 Assume $E$ is the symbol to the parametrix to the operator corresponding to $f$,
 that is $Ef \sim I$ (modulo regularizing action), then $\mbox{ ker }E$ is a polar set.
 \end{parametrix}

 Define $(I_{1})$ such that $U^{\bot} \varphi=0$ for 
 $\varphi \in (I_{1})$, that is a polar. Assume $U_{1}^{\bot}(p \varphi)=U^{\bot}(\varphi)=0$, 
 given $\varphi$ hypoelliptic and $p$ algebraic, we have that $U^{\bot} \varphi=0$ implies $\varphi=0$. 
 Conversely, the support of $\varphi \in \mathcal{D}_{L^{1}}$ is a translation domain. Assume
 that we have existence of $E$ such that $E(\varphi)=I$, when $\varphi \neq 0$
 and $E(0)=0$. If we assume $U$ surjective on a symbol ideal $(I)$ and $E U^{\bot} = U^{\bot} E$,
 in particular $U I = I U$, we have that $(I_{1}) \sim \mbox{ ker }E$, that is the kernel to 
 the parametrix can be identified with a polar.
 
  Assume $Uf \in \mathcal{D}_{L^{1}}$ and $U^{*} \widehat{f}$ analytic. The graph norm 
  to $\mathcal{G} \times \mathcal{G}^{\bot}$ is given by $\parallel U f \parallel_{1} + \parallel U^{*} \widehat{f} \parallel_{1} \simeq \parallel Uf \parallel_{G \times G^{\bot}}$.
  Thus, given the equation above, we have that sets invariant for graph norm implies presence of
  polar sets, which implies presence of sets invariant for $U$ In $L^{1}$. When 
  $U^{*} \widehat{f}=0$ on a set $\Sigma^{*}$, the equality can be continued to $L^{1}$ 
  using density and $\Sigma^{*}$ can be determined using continuity. Further, $\Sigma^{*} \rightarrow \Sigma$ 
  can be determined using duality (Fourier). Note that if $U f=U^{\bot} f=0$
  on $\Sigma$, we have that $f=0$ on $\Sigma$. Given $d U=a d x + b d y$ and 
  $d U^{\bot}=a_{1} d x + b_{1} d y$ and $\phi : a / b \rightarrow - \overline{b_{1}} / \overline{a_{1}}$ 
  projective, we have that $d U \bot d U^{\bot}$ in the sense of differential forms 
  (wrt $\int d U d \overline{U^{\diamondsuit}}$ \cite{AhlforsSario60}). Given $U$ harmonic we have 
  $\Sigma^{*} \sim \phi(\Sigma)$. When ${}^{t} U \rightarrow U^{\bot}$ is completed to $L^{1}$ and harmonicity 
  by an algebraic continuation ($f$ HE), we thus have that $U$ is projective for $f$.
  Assume $\frac{d U_{2}}{d U_{1}}=\alpha_{2}$.
 
 \newtheorem{generator}[adj]{Lemma Transversal generator}
 \begin{generator}
 Given $\Omega$ a $U_{1}-$ domain and $g / \alpha_{2} \in L^{1} (d U_{2})$, then $\Omega$ 
 is a $U_{2}-$ domain relative $g/ \alpha_{2}$. Given a symbol ideal $(I)$, such that 
 $(I) / \alpha_{2} \subset (I)$, we have that every $U_{1}-$ domain for $(I)$ is simultaneously a 
 $U_{2}-$ domain. In $L^{1}$ we have that every (closed) domain $\Omega$, can be given as a $U_{1}$ 
 domain relative a generator $f_{0} \in L^{1}$ locally (\cite{Garding98}, Chapter 6). 
 \end{generator}
 
 Absence of essential singularities in $\infty$ (finite order) means that every analytic 
 function monotropic to $ f_{0} \in L^{1}(d \lambda)$ can be represented as $g / \alpha$, 
 where $g,\alpha$ are entire (\cite{Cousin95}). Note that
 $\{ \tau \varphi \}$ is not relatively compact in $ \mathcal{D}_{L^{p}}$ (\cite{Schwartz66}). Thus when 
 $\{ \varphi < \lambda \}$ relatively compact, this does not imply that 
 $\{ \tau \varphi < \lambda \}$ is relatively compact. 
 A sufficient condition for $\subset \subset$ to be preserved, is that the movement is 
 downward bounded. When $(I)=(I_{HE})$, we can choose $U$ 
 as monotonous sequential movements in one parameter.

 Assume $V$ a movement analytic over $\gamma$ and $U$ a continuous movement, such that 
 $\mid {}^{t} (U-V) \gamma \mid \leq \epsilon$. Assume $U_{1}$ a harmonic movement and 
 $d V / d U_{1}=e^{\phi}$. Starting from $V$, we can regard $U$ as a continuous monotropic 
 continuation of $V$. Where $V$ changes into $U$, we assume $\phi$ linear in $x,y$. 
 The condition for monotropy is thus, that we can always locally find an analytic movement, 
 with $U \gamma \sim V \gamma$.  Assume the boundary corresponds to $C$ and given 
 $V$ analytic, we can assume $d V$ closed. Given $V^{\blacktriangle}$ reduced, we have that $d V$ 
 is exact. Thus an analytic sequential movement can be continued using a harmonic movement, 
 to a closed movement. 
  
 \subsubsection{The orthogonal distribution}

 Assume the polar semi algebraic, with compact
 sub level surfaces for $\nu + i \vartheta$, then we have that the sets in $\nu,\vartheta$ 
 are semi algebraic, but do not both have compact sub level surfaces. 
 When we consider $U$ as distribution, with support in the polar $\Omega$, there is a movement $V$ 
 with support in the $\complement \Omega$.
 
 Given that we have an algebraic base for the topology, we can have absence of
 trace (as in $\tilde{C}$) with possible presence of lineality. Otherwise, this can not be determined.
 
 For phe operators, we have that the complement to the range to the symbol, can be represented on 
 a spiral domain (\cite{Dahn18}).
 It is a necessary condition for hypoelliptic symbols that the ``polar'' does not contain a spiral.
 Further, $G \times G^{\bot} (I) \subset (I)$ with $(I)=(I_{HE})$ implies $G \neq G^{\bot}$. 
 The condition for the polar is written $G \times G^{\bot} (I)^{\bot} \subset (I)^{\bot}$ 
 implies $G \neq G^{\bot}$. In particular 
 we assume the polar generated by $G \times G^{\bot}$ (euclidian). Example: when $G=G_{1}$
 and the polar is a hyperboloid, the polar has zero dimension!

 \newtheorem{ortogonal}[adj]{Lemma An orthogonal movement}
 \begin{ortogonal}
 Consider $d w \bot d (U,U^{\bot})$ implies $d w=0$ on a set $K$ (defined by $\tilde{C}$). 
 Thus, given $w$ absolute continuous and $w^{\bot}=I-w$ we have that $d w=0$ implies that $w$ has support on $K$,
 that is defines an orthogonal distribution.
 \end{ortogonal} 
 \subsection{The boundary}

 When $U$ is analytic and $\gamma$ in an annulus, we can define $\nu \rightarrow \vartheta$ 
 as projective. When $U$ is harmonic the mapping can be continuously continued to a disk. Assume for this reason 
 $d U / d U_{1}=\alpha$ real, then we have that where $\alpha \neq 0$, the domain can be considered 
 as contractible. When we assume $d U_{1}=0$ implies $U_{1}^{\blacktriangle}=0$, that is isolated points, 
 harmonicity can be continued over isolated points and the domain for $U_{1}$ can be assumed 
 symmetric.

 \newtheorem{higher}[adj]{Lemma Movements of higher order}
 \begin{higher}
 Assume the boundary is given by $C$, points where the movement changes character. 
 Over $C$, we have that $U^{\blacktriangle}=0$, that is we assume the movement changes character 
 through $I$. Assume $\Sigma = \{ X(f)=0 \}$, that is a domain where $U$ is analytic. 
 The condition $\exists U_{k} \quad \Sigma^{(k)} \neq \{ 0 \}$ means that there is a 
 movement, analytic over $f$.   
\end{higher} 

 Assume $L^{1} \ni f \rightarrow Uf \in L^{1}$, where we assume $U \in \mathcal{D}_{L^{1}}'$ 
 that is with algebraic base. Given a very regular boundary, it is to determine the character of the 
 movement, sufficient to consider the movement in the phase space and we can assume 
 $f=e^{\phi}$ with $\phi \in L^{1}$ and $\log U f \in L^{1}$. We assume
 existence of movement, possibly of higher order, analytic close to the boundary. 
 
 \newtheorem{harmconj}[adj]{Proposition Harmonic conjugation of movement}
 \begin{harmconj}
 Note that given $X(f)$ to $U^{\diamondsuit}$, we have that 
 $U^{\diamondsuit}$ is harmonic if $X(f)=\delta X(f)=0$. Then $U \rightarrow {}^{t} U$ 
 preserves character of the movement, given that $X$ has coefficients algebraic in $x,y$. 
 This means that the derivatives for $U^{\diamondsuit}$, that is $\xi,\eta$ are real 
 and analytic. Given that some derivative for f is non zero (cf, envelop), the singularities 
 to $U^{\diamondsuit} f$ are zero's to $\xi,\eta$. Note that given 
 $d V \bot d U^{\diamondsuit}$ and that $d U$ is closed, we have that $d V$ is exact, 
 assuming $(U^{\diamondsuit}, d f)=\int f d U^{\diamondsuit}$ and that $f$ has compact support.
 \end{harmconj}

 (\cite{AhlforsSario60})
 In particular $U^{\bot} \bot d f$ implies $(U^{\bot}, d f)=\int f d V^{\bot}=0$, 
 for a movement $V$. Given $\xi,\eta$ algebraic and $w$ a differential form, 
 $(w, d U^{\diamondsuit} f)=\int U^{\diamondsuit} f d w=(w , (-\eta + i \xi) d f)=
 ((-\eta + i \xi) w, d f)$.

 Note that if $f$ analytic and in $\mathcal{D}_{L^{2}}$, we have that $f \sim (U + U^{\bot})f$, 
 when $Uf$ is analytic and $U \in \mathcal{G}$. 
 We assume $U^{\bot}$ the completion in $\mathcal{D}_{L^{2}}$, such that $<U^{\bot} f, U f>=0$.
 Using the projection theorem, $f \sim U f + U^{\bot} f_{1}$ in $L^{2}$. Let $U^{\bot} \rightarrow I$
 regularly, when $U^{\bot} f_{1}$ non-trivial.
 Assume now that $U$ is not projective,
 that is consider $U f + U^{\bot} \widehat{f}$, when $\widehat{f} \in \dot{B}$. 
 
 Given $f \in \mathcal{D}_{L^{1}}$
 we assume $\widehat{f} \in (\dot{B})$, why we can consider $U^{\bot} I \in \mathcal{D}_{L^{1}}'$
 (given nuclearity a measure).  We use that 
 $d U^{\bot}/ dU=p$ is a polynomial locally. 
 Consider $<f,\widehat{g}>$ for $f,g \in \mathcal{D}_{L^{1}}$. Assume $U_{1}$ harmonic translation
 with $d U_{1}^{*} \rightarrow 0$ implies $d \lambda \rightarrow 0$ (Lebesgue). Assume there is
 a $V$ with $d V / d U_{1}^{*}=\alpha^{*}$, where $\{ \alpha^{*} < \mu \} \subset \subset \Omega$
 (globally), $\mu$ constant. Assume $f \bot \widehat{g}$, with $f,g$ non trivial. Thus, if $V \rightarrow 0$ regularly, we have $U_{1} \rightarrow I$
 regularly. That is $N(V)=R(U_{1})^{\bot}$ (zerospace and range). Note that if $d U_{1}^{*}/ d \lambda=\lambda_{1}$ and $d U_{1}^{*} / d U_{1}=\rho$,
 given $\rho / \lambda_{1} \neq constant$, then a sufficient condition for a regular approximation is 
 $\alpha^{*} \neq constant$. Consider now $<(U + U^{*} + V) f, \widehat{g}>=<f,{}^{t} (U + U^{*}) \widehat{g}>$.
 Then, when ${}^{t} U$ $(\neq {}^{t} U^{*}$) is projective over $\widehat{g}$, we can continue
 $U^{*}$ with ${}^{t} V \widehat{g} \rightarrow 0$. When $V$ is assumed reflexive, we have
 that $V \in \mathcal{G}^{\bot}$ and $V$ defines an orthogonal distribution with support on $\mbox{ ker }{}^{t} V$.
 When ${}^{t} U^{*} \widehat{g}=0$ over a set (polar), we extend the kernel to ${}^{t} V$ with
 this set. In this manner $U^{\bot}$ can be defined in $(\mathcal{D}_{L^{1}})'$.
 Note that maximal rank does not imply that the movement can be uniquely determined. 
 Starting from $U + U^{\bot}=I$, we have that $U_{1} + U_{2}$ have the same (maximal) rank 
 that $U_{2} + U_{1}$, that is the condition on constant rank does not identify the movement. 
 
 Assume the defect spaces are defined by $D(V)^{\bot}$ and $R(V)^{\bot}$, where $V \sim \log S$. 
 A closed symmetric extension is maximal iff the defect spaces have the same dimension. 
 Given $S \sim (U,U^{\bot})$ symmetric, that is $U$ reflexive and densely defined, we consider 
 the continuation to $(U,U^{\triangle})$. Here we consider $(U,U^{\bot})$ as a sequential 
 movement, that is a 1-form and $(U^{\bot},U) \simeq (U,U^{\triangle})$ is
 seen as a ``interpolations property''. When the continuation is algebraic, that is 
 corresponding to a removable set, we have a global base for the movement, given the range is dense. 
 The condition $(I U) =(U I)$ means that $U^{\bot} \gamma=0$ iff 
 $U^{\bot}(- \gamma)=0$. Algebraicity implies the defect indexes are equal.
 Symplecticity implies the defect indexes zero.

 Assume $\mathcal{G}$ is constructed using $q_{j}$, quasi orthogonal, with distinct zero's. 
 With these conditions, we can localize domains for absolute continuity for a measure. 
 Assume for this reason $d U_{j}=q_{j} d U_{1}$, where $q_{j}$ have distinct zero's,
 forming a base for absolute continuity. When for instance $d U_{3}=p dU_{2}$, we have that 
 $p = q_{3}/q_{2} \neq 0$, where $q_{2},q_{3}$ have no zero's in common. If $\int f d U_{j}=0$ 
 we have that $\int_{\Omega} f q_{j} d x=0$, why $f q_{j}=0$ on $\Omega$. Given $q_{j} \neq 0$ 
 on $\Sigma \cap \Omega$, we must have $f=0$ on $\Sigma \cap \Omega$.
 Note that we can choose $q_{j}$ orthogonal, if $f \neq 0$ on $\Sigma \cap \Omega$. (\cite{Riesz23})
 
  \subsection{The boundary to the completed movement}

According to (\cite{Lie91}, Chapter 6), we have that $A f \equiv 0$ determines the movement $X(f)$ iff $X(Af)-A(X(f)) = \lambda A f$.
For instance rotation, translation have $(X A) \equiv 0$. Given $A f \rightarrow A^{*} \widehat{f}$
we have that $A^{*} \widehat{f} \equiv 0$ determines a movement. Assume $X^{\triangle} A^{\triangle}f  - A^{\triangle} X^{\triangle} f = \lambda A^{\triangle} f$.
Given $\frac{\eta^{\triangle}}{\xi^{\triangle}} \sim \frac{\eta^{*}}{\xi^{*}}$, we have that consequently 
$(X^{\triangle} A^{*}) \equiv 0$. Assume $X^{\triangle}$ symplectic and $X^{\triangle}=(X_{1} X^{*})$, if $(X_{1} A^{*}) \equiv 0$ and $X_{1}(0)=0$
we have that $(X^{\triangle} A^{*}) \equiv 0$.
Note that when $(\xi,\eta) \sim (\alpha,\beta)$ (translation) gives $\eta' \equiv 0$
for the continued group (see last section) and when $(\xi,\eta) \sim (-y,x)$, we have that $\eta' \equiv (1 + (y')^{2})$.

Assume $\chi^{*} \rightarrow \chi^{\triangle}$ algebraic and $\mathcal{G}^{\bot} \simeq \mathcal{G}^{*}$,
we thus assume $\mathcal{G} \times \mathcal{G}^{\bot} \rightarrow \mathcal{G} \times \mathcal{G}^{\triangle}$
algebraic, that is $\mathcal{G}^{\triangle} \simeq e^{P(x^{*},y^{*})} \mathcal{G}$. Further,  given
hypoellipticity, $\mathcal{G},\mathcal{G}^{\bot}$ can be defined separately, when $\mathcal{G} \neq \mathcal{G}^{\bot}$.

Consider $\mathcal{T} X^{\bot}(f)=
\widehat{X} \frac{\delta \widehat{f}}{\delta x^{*}} + \widehat{Y} \frac{\delta \widehat{f}}{\delta y^{*}}$.
Assume $\frac{\delta}{\delta x}(\widehat{U f})=\xi^{*}$ and $\frac{\delta}{\delta x} U^{\bot} \widehat{f}=\widehat{X}$.
Given $(\xi^{*},\eta^{*}) \bot (-\widehat{Y},\widehat{X})$ and $\xi^{*} / \eta^{*} \sim \widehat{X} / \widehat{Y}$,
$U$ can be determined from $U^{\bot}$. Given $U^{\bot}=(I-U)$, we have that $X^{\bot}(f)=X_{1}(f)- X(f)$,
that is when $X^{\bot}=0$ we have that $U$ corresponds to translation.

 Assume
 $M=x \frac{d}{d x}$, then we have that $-\mathcal{L}(M f)=\widehat{f}(x^{*}) +  M \widehat{f}(x^{*})$. Thus, if
 $\mathcal{T} M =- M \mathcal{L}$, we have given $M^{\bot}= I - M$, $\mathcal{L} (-M)^{\bot} = \mathcal{T} M$.
 Note that $(M f)^{\bot}$ can be defined by $<f, \phi + x \frac{d \phi}{d x}>=0$. Sufficient for this
 is that $f \bot \phi$,$f \bot \frac{d \phi}{d x}$. 

\newtheorem{rambd}[adj]{Lemma Completion by ramified boundary}
 \begin{rambd}
Assume $C$ discrete and $\psi(C)=S$ connected, that is $S$ corresponds to a covering $\tilde{C}$.
For instance $d U_{2} / d U_{1} = \alpha_{2}=const$ on $S$. Let $\mid d U \mid  \sim d U(\nu,\vartheta)$.
When we consider $d U_{1} \rightarrow d U_{1}^{\diamondsuit}$, this corresponds $\nu \rightarrow - \vartheta$.
$S$ for this reason describes a set of symmetry for $d U(u,u^{\diamondsuit})$. 
Further, given $U^{\triangle}=U$
describes the completion of $U^{*}$ in $L^{1}$, we assume the completion maps $C \rightarrow \tilde{C}$.
\end{rambd}

Assume the movement $V$ defined by $A F \equiv 0$ and $\{ {}^{t} V,F \}=A F$.
The condition $X_{j}(A f)=\lambda A f$, can be written $X_{j}(X) f_{x} + X_{j}(Y) f_{y} + X X_{j}(f_{x}) + Y X_{j}(f_{y})=\lambda (X f_{x} + Y f_{y})$.
A sufficient condition for this is that $\frac{X_{j}(X)}{X} \sim \frac{X_{j}(f_{y})}{f_{y}}$ and $\frac{X_{j}(Y)}{Y} \sim \frac{X_{j}(f_{x})}{f_{x}}$.
Given $A f \equiv 0$ on a domain, where $f$ is analytic, we have that $\frac{d y}{d x} \sim -\frac{\delta f / \delta x}{\delta f / \delta y}$.

 \section{Sums of squares}
  Assume $U_{1}$ translation and $\frac{d U_{j}}{d U_{1}}=\alpha_{j}$ regular for all j 
 (without constant surfaces). Assume further $\frac{d U_{1}}{d \lambda(x)}=\frac{1}{\beta}$ regular 
 for the Lebesgue measure. Then we have that $\int d U_{1} \ldots d U_{r}=\int \alpha_{2} \ldots \alpha_{r} d U_{1}= \int (A(X)/ \beta) d \lambda (x)$.
 Given $A/\beta$ polynomial, we have that $\int_{V} (A/\beta) d \lambda (x)=0$ implies V of measure zero (Hurwitz).

 \subsection{Topology}
  Assume $f(\zeta_{T})=U f(\zeta)$, given isolated singularities and Schwartz type topology, 
  where $\zeta_{T}$ is translation domain, we have that $U$ is translation locally.  
  The lineality $\phi(\zeta_{T}) \equiv \phi(\zeta)$ then corresponds to $U_{1} \phi \equiv \phi$. 
  The corresponding movement on $f=e^{\phi}$ has the same character. The condition $U_{1}^{\bot} \equiv 0$ 
  iff $U_{2}^{\bot} \equiv 0$, is satisfied if $d U_{1}^{\bot} /d U_{2}^{\bot}=konst.$
  Thus, given $\alpha_{2} d U_{1} = d U_{2}$ and when the lineality is separated from the trace $\tilde{C}$
  (assume $d U_{1} \neq 0$), we have that $\alpha_{2}$ is not constant on sets of positive measure.

  Consider $U F(\gamma)(\zeta) = F({}^{t} U \gamma)(\zeta)$ $\rightarrow {}^{t} U \gamma$ 
  dependent of $(x,y,\frac{d y}{d x})$. Further, consider ${}^{t} U \gamma \rightarrow {}^{t} U_{1} \gamma$. 
  Through the condition on $\alpha=e^{\phi}$, we can assume ${}^{t} U_{1} \gamma \rightarrow {}^{t} U_{0} \gamma$, 
  where $U_{0}$ is harmonic translation, analytic and real. Assume $\psi$ a proper mapping 
  $C \rightarrow \tilde{C}$ that generates a covering of $C$. Thus, $U_{j}^{\blacktriangle}=0$ can 
  be related to $U_{0}^{\blacktriangle}=0$. 
 
  Define $\Sigma_{0}=\{ \delta X(f)=0 \}$ and $\Sigma=\{ X(f)=0 \}$, thus given $X$ defines
  $U^{\bot}$, we have that the edge of the envelop is given by harmonic movements. Given 
  $U^{\bot} \rightarrow U^{\triangle}$ to harmonicity, we have that the converse mapping 
  is multivalent. Assume $\pi : \Sigma \rightarrow \Sigma_{0}$, given a set $\Omega \cap \Sigma = \emptyset$, 
  we still have that  $\tilde{\Sigma_{0}} \cap \Omega \neq \emptyset$,
  thus $U^{\triangle}$ does not determine the movement $U$ uniquely (\cite{Julia55}), but there is always  
  some $V \in \mathcal{G}$ not necessarily harmonic, such that $V^{\triangle}=U^{\triangle}$. 
  Consider 
  $U \rightarrow U^{*} \rightarrow U^{\triangle} \rightarrow U^{\bot}$, then we have existence 
  of $V$ such that $U^{\bot} \sim I - V = V^{\blacktriangle}$.
 
  Starting with $\mathcal{G}(I) \subset (I)$ and $d U = \alpha d U_{1}$, we have that a 
  necessary (and sufficient) condition for inclusion between weighted ideals, is that 
  $1/\alpha \rightarrow 0$ in $\zeta - \infty$. We assume $\Sigma$ the set of $(x,y)$, where 
  $\alpha$ is defined locally. $\Sigma$ is divided by $C$, the set of $(x,y)$ where the 
  movement changes character. Assume $\alpha=\alpha_{j}$ on $\Sigma_{j}$ and 
  $\Sigma = \cup \Sigma_{j}$ and denote $\tilde{\alpha}_{j}$ for the continuation to 
  $x,y - \infty$. Note that $\frac{1}{\alpha + 1/ \alpha} \rightarrow 0$,
  when $\alpha \rightarrow 0$ and when $\alpha \rightarrow \infty$.  In 
  particular we have that $\Sigma_{j} \alpha_{j} + 1/\alpha_{j} \rightarrow \alpha$ when 
  $\mid \zeta \mid \rightarrow \infty$. Assume $(\alpha,\beta)$ defines
  a global base for the measures and that $U \rightarrow U^{\bot}$ projective, then $(\alpha,\beta)$ defines inclusion between the ideals 
  globally. The choice of $\beta$ is dependent on conjugation and for this reason the choice 
  of norm in the condition for inclusion.  
 
 \subsection{Projectivity for movements}

 Assume $X(f)$ defines $Uf$ and $X^{*}(\widehat{f})$ defines $U^{*} \widehat{f}$,
 then we have that $\mathcal{T}(X_{f})=\xi^{*} \frac{\delta \widehat{f}}{\delta x^{*}} + \eta^{*} \frac{\delta \widehat{f}}{\delta y^{*}}$,
 where $\mathcal{T}$ completes $x \rightarrow x^{*}$ to symplecticity. Assume for this reason 
 $U^{**} \rightarrow U$ locally injective, when the propagation direction is fixed, we have 
 that $\mathcal{G}$ nuclear, that is given $\overline{\delta}(Vf)=Uf=0$, we have that 
 $V f=\int f d \mu$ where $d \mu=d \mu_{1} \otimes \ldots \otimes d \mu_{r}$
 and $V_{j} \in \mathcal{G}$.

 When $\Sigma$ is defined by the principal part of the symbol, $\Sigma=\{ X_{j}=0 \}$
 is given by lines (planar curves are lines).
 As above, we can map $X_{j}=0$ on $X_{j}^{*}=0$ on $X_{j}^{**}=0$, where lines are mapped 
 on to lines. In the plane, we have that if $X_{j}=0$ can be given by a second order d.o., 
 the trajectories are zero lines.

  Assume $X,Y$ have compact sub level surfaces (analytic on pseudo convex domains). Assume 
  $\phi(\frac{d y}{d x})=\frac{y}{x}$. 
  Assume $Y/X \rightarrow X/Y$ projective. 
  Using Radon Nikodym (\cite{Riesz56}), if $d M(\frac{y}{x})=N \phi^{-1}(\frac{x}{y})$ and 
  $N(\frac{d x}{d y})=0$ implies $d M (\frac{y}{x})=0$, we have that 
  $d M (\frac{y}{x})=N(g \frac{d x}{d y})$, for some $g$ measurable with respect to $N$. Given $d M \sim d M^{\diamondsuit}$, 
  we have that $d M$ harmonic (real, \cite{AhlforsSario60}), that is with zero's such that where
  $\mid y / x \mid = 1$, we have $d M(\nu,\vartheta) \sim N(g \frac{d y}{d x}, g' \frac{d x}{d y})$.
  
 \newtheorem{minimal}[adj]{Proposition Minimal characteristics}
 \begin{minimal}
 Assume $\overline{\delta}(V f)=X(f)=0$ defines $\Sigma$
 and that $\Sigma_{0}$ is defined by $\delta X(f)=0$. Then we have that 
 $\Sigma_{0}$ is minimal if $V f$ harmonic. Given $p$ a planar 
 geodetic in $\Sigma_{0}$, such that $d p^{\bot}= d \gamma$, where $\gamma$ 
 is assumed algebraic, we have that $\Sigma_{0}$ is algebraic. Given 
 $(x,x^{\triangle}) \rightarrow (x^{\triangle},x)$ isometrically isomorphic (projective
 completion to a normal model) and consider $p = (x,x^{*}) \mid_{L}$ (restriction to a plane)
 such that $d p^{\bot}=d \gamma$, where $\gamma$ algebraic, then we have that a minimal domain
 for continuation $x^{*} \rightarrow x^{\triangle}$, is algebraic.
 \end{minimal}
  (Cf. \cite{Lie78})

\newtheorem{pureortho}[adj]{Lemma The graph has a pure orthogonal}
\begin{pureortho}
 When we consider $\mathcal{G} \times \mathcal{G}^{\bot} \backslash (\mathcal{G} = \mathcal{G}^{\bot})$,
 we can consider $g + i g^{\bot}$. Assume reflexivity such that $g^{\bot \bot} \simeq -g$, 
 then we have that $(g + i g^{\bot})^{\bot}=-i(g + i g^{\bot})$. 
 \end{pureortho}
 
 Assume for instance $g_{1}^{\bot}=(g_{2} + \ldots + g_{r})$.
 Assume $d g_{2} / d g_{1}=\alpha_{2}$ with compact sub level surfaces (regular). 
 Assume $g^{\bot} \rightarrow g$ locally 1-1, then we can identify movements on an interval 
 in $\mathcal{G}$. Assume $g \rightarrow g^{\bot}$ through a cylindrical domain, 
 then $g^{\bot}$ can be identified as a movement in $\mathcal{G}$, 
 where $d g^{\bot}/d g_{1} \neq 0$.

 Assume $\overline{\delta}(U^{\diamondsuit} f)=U^{\diamondsuit}X(f)$ and $<\overline{\delta} (U^{\diamondsuit}f),\psi>$
 $= <U^{\diamondsuit}f, \overline{\delta} \psi>$, that is given $U^{\diamondsuit} f$ analytic, we have 
 that $U^{\diamondsuit} f \bot \overline{\delta} \psi$. In the same manner, given $< \xi f, \psi>=<f, \xi \psi>$, 
 we have $<X(f),\psi>=<f , X(\psi)> - <f, (\frac{\delta \eta}{\delta x} - \frac{\delta \xi}{\delta y}) \psi>$
 
 \emph{Thus the character for $U^{\diamondsuit} \rightarrow {}^{t} U^{\diamondsuit}$ is preserved, given that 
 $f \bot (\frac{\delta \eta}{\delta x} - \frac{\delta \xi}{\delta y}) \psi$.
 When $\frac{d y}{d x} \sim \frac{d \eta}{d \xi}$, we have formally 
 $\Delta U^{\diamondsuit}=0$.}
 
 Assume $<U^{-1} \frac{\delta}{\delta x} U f,g>=< U f, \frac{\delta}{\delta x} {}^{t} U^{-1} g>$.
 Let ${}^{t} U \frac{\delta}{\delta x} {}^{t} U^{-1} g=\xi$. Thus $<X_{U}(f),g>=<f,X_{{}^{t} U^{-1}}(g)>$.
 Further, $<\xi \frac{\delta f}{\delta x} + \eta \frac{\delta f}{\delta y},g>=
 <f, \xi_{1} \frac{\delta g}{\delta x} + \eta_{1} \frac{\delta g}{\delta y}>$. Assume 
 $<\xi f,g>=<f, {}^{t} \xi g>$. Thus, $<X_{U}(f),g>=
 <f, (\frac{\delta {}^{t} \xi}{\delta x} + \frac{\delta {}^{t} \eta}{\delta y}) g > + 
 <{}^{t} \xi \frac{\delta g}{\delta x} + {}^{t} \eta \frac{\delta g}{\delta y}>$.
 Thus given the compatibility condition, $\frac{\delta {}^{t} \xi}{\delta x} + \frac{\delta {}^{t} \eta}{\delta y}=0$
 (mass conservation), we have that ${}^{t} \xi = \xi_{1}$,${}^{t} \eta = \eta_{1}$.

 Given $d U^{\bot}$ locally algebraic, we have that its zero set is removable (cf. strict pseudo 
 convexity). In the case with $r_{T}'$ (\cite{Dahn13}), we use $F(r_{T}' \gamma)=f(\zeta_{T})$, that is 
 $U$ is acting on $x,y$. The condition on strict pseudo convexity, is $\Delta_{L} F > 0$ 
 and $\alpha=\frac{\delta x}{\delta \zeta_{1}}$ and $\beta = \frac{\delta y}{\delta \zeta_{1}}$, 
 where $\zeta_{1},\zeta_{2}$ generates a complex line $L$ and where we assume $\alpha,\beta$ 
 real and not both zero. The condition on $\alpha,\beta$ is interpreted as 
 $\gamma$ regular on $L$. (\cite{Oka60})
 
 \newtheorem{Thorin}[adj]{Proposition Maximal subspace for projectivity}
 \begin{Thorin}
 There is a subspace $\mathcal{D}_{L^{r}} \subset \mathcal{D}_{L^{2}}$,
 such that $U \rightarrow U^{\triangle}$ is projective. 
 \end{Thorin}
 
 Assume $f \in L^{1}(d U)$ and $\widehat{f} \in L^{2}(d U)$ and that $U U^{\triangle}=U^{\triangle}U$
 in $L^{2}$ (a normal operator). Assume $T(f)$ is the continuation
 of $f$ according to $\int \tilde{f} d U \sim \int f d \tilde{U} \sim \int f d U^{\triangle}$ 
 (cf. the continued group). Then we have that $\parallel T(f) \parallel_{r} \leq 
 \parallel T(f) \parallel_{1}^{1-\theta} \parallel T(f) \parallel_{2}^{\theta}$,
 where $0 < \theta < 1$ (cf Malgrange (1957)). 
 Given the movement is reflexive and where it is absolute continuous,the inequality holds for both $U,U^{\triangle}$.
  A polar set in $L^{2}$ corresponds to a polar set in $L^{r}$. 
  Assume $U \rightarrow U^{\triangle}$ projective relative $L^{2}-$ norm, then in $L^{r}$,
 $\parallel f \parallel_{r} \sim_{0} \parallel Uf \parallel_{r} + \parallel U^{\triangle} f \parallel_{r}$
 ($\sim_{0}$ means geometric equivalence), where $U \neq U^{\triangle}$.
 In particular,
 when $U^{\triangle} f=0$ in $L^{1}$ and using Parseval in $L^{2}$,
 we have that $U^{\triangle} \sim_{0} 0$ in $L^{r}$ implies $U \sim_{0} I$ in $L^{r}$
 
 Note that when $f$ analytic and $d U \rightarrow d U^{\triangle}$ projective (the measures
 are absolute continuous respectively), then $<d U,f>+<d U^{\triangle},f>=0$. When $<U,d f>=0$,
 given $<I, d f>=f$, that is $U f$ absolute continuous, when $U \rightarrow I$, we have
 $<U^{\triangle} + U, d f>=<I, d f> \sim f$. We can write $L^{1}_{ac}(d U) \subset L^{r}(d U)$.
 Conversely, when the movement is harmonic (or $\delta X_{U}(f) \geq 0$), why $U \rightarrow {}^{t} U$
 preserves character, we can derive absolute continuity for $U$. When $f$ hypoelliptic with
 $E$ the symbol to the parametrix, we have $E({}^{t}(U + U^{\triangle} - I) f)=0$ and since
 the kernel to $E$ is trivial, ${}^{t} U$ is projective over $f$.
 \subsection{Spiral domains}
 
 Assume $U_{1}$ translation and $d U_{s}=\alpha_{s} d U_{1}$. A continuation according to 
 $U_{s}$ corresponds to a non-closed extension, that is it does not preserve algebraicity. 
 For a closed extension, it is necessary that the Dirichlet integral is finite. For $d U=0$ 
 to imply Lebesgue measure zero, it is sufficient that algebraicity is preserved.
 Thus, $\alpha_{s}$ is continuous, but not algebraic. 
 
 \newtheorem{hecond}[adj]{Lemma Change of sign in infinity}
 \begin{hecond}
 $\mbox{ Re }P$ changes sign on a component $\ni \infty$ implies $\Delta(P) \neq \{ 0 \}$.
 \end{hecond}
 (\cite{Dahn15}, Lemma 2.1)
 
 For the spiral, we have that $\xi_{s}=(y + \kappa x)$,$\eta_{s}=(-x + \kappa y)$ 
 and the translation coefficients, $\xi_{1}=c_{1}$,$\eta_{1}=c_{2}$ are constants. Given $\frac{d U_{s}}{d t}=0$, we have that 
 $\eta_{s} / \xi_{s} \sim const$. Obviously $\xi_{s},\eta_{s}$ are 
 polynomials in $x,y$. Finally, $\frac{d U}{d U_{1}}=\frac{(y + \kappa x) + (-x + \kappa y) \frac{d y}{d x}}{c_{1} + c_{2} \frac{d y}{d x}}$
 and $\sim_{0} (-x + \kappa y) \frac{d x}{d y} + (y + \kappa x)$.

  \newtheorem{spiraldiag}[adj]{Lemma The spiral transformation is on the diagonal}
 \begin{spiraldiag}
 Assume $\overline{X}^{\diamondsuit}(f)=0= (-x+k y) \frac{\delta f}{\delta x} + (y + k x) \frac{\delta f}{\delta y}$. 
 When we have that $-\mathcal{L}\overline{X}^{\diamondsuit}(f)=(-x^{*} + k y^{*}) \frac{\delta \widehat{f}}{\delta x^{*}} + (y^{*} + k x^{*}) \frac{\delta \widehat{f}}{\delta y^{*}}$
 Thus when $\overline{g_{s}}^{\diamondsuit} \in \mathcal{G}$, we have that ${\overline{g_{s}}^{\diamondsuit}}^{*} \in \mathcal{G}$ 
 and the spiral is on the diagonal $\overline{g_{s}}^{\diamondsuit} \sim {\overline{g_{s}}^{\diamondsuit}}^{*}$
 in $\mathcal{G} \times \mathcal{G}^{\bot}$. Assume $V f \equiv (y + \kappa x) \frac{\delta f}{\delta x} + (-x + \kappa y) \frac{\delta f}{\delta y}$
 (\cite{Lie91}) when $\widehat{(y \frac{\delta f}{\delta x})}=-x^{*}\frac{\delta \widehat{f}}{\delta y^{*}}$
 and so on, we have that $-\widehat{V f}= 2 \kappa \widehat{f} + 
 (-y^{*} + \kappa x^{*}) \frac{\delta \widehat{f}}{\delta x^{*}} + (x^{*} + \kappa y^{*}) \frac{\delta \widehat{f}}{\delta y^{*}}$
 Thus, except for the first term (modulo a translation) also $\widehat{Vf}$ defines a spiral (with $-\kappa$).
 
\end{spiraldiag}
 
 Note, if $X_{s}$ denotes the spiral $U_{s}$ as above, we have that $d U_{s} \simeq (X_{s} + \kappa) d t$.
 Assume $f$ Hamiltonian to $X,Y$, that is $A(g)= \{ f,g \}$ and $h=A(g)$. Then we have that 
 over an involutive set, that is $\{ f,g \}=0$ that $X_{s}(h)=0$ iff $X_{s}(h) + \kappa(h)=0$.
 For the spiral, we have that $V_{s} \simeq V_{s}^{*}$, that is $X_{s} \simeq X_{s}^{*}$.  
 
 Consider $X_{s}(X_{s}(f))=\kappa X_{s}(f) - X^{\diamondsuit}_{s}(f) + \xi_{s} X_{s}(\frac{\delta f}{\delta x}) + \eta_{s} X_{s}(\frac{\delta f}{\delta y})$.
 Assume $g=U_{s}f$, then $\xi X_{s}(\frac{\delta g}{\delta x}) + \eta_{s}(\frac{\delta g}{\delta y})=
 (\xi_{s}^{2} + \eta_{s}^{2}) X_{s}(f)$ and $\kappa X_{s}(g) + X_{s}^{\diamondsuit}(g)=
 \kappa (\xi_{s}^{2} \frac{\delta f}{\delta x} + \eta_{s}^{2} \frac{\delta f}{\delta y}) + 
 - \xi_{s} \eta_{s} \frac{\delta f}{\delta x} + \xi_{s} \eta_{s} \frac{\delta f}{\delta y}$.
 Thus, when $X_{s}(f)=0$ and $\frac{\delta f / \delta x}{\delta f / \delta y}=\frac{\eta_{s}}{\xi_{s}}$,
 then $X^{2}_{s}(f)=0$ iff $\frac{\xi_{s}}{\eta_{s}}=-\frac{(k-1)}{(k+1)}$. Obviously, $X_{s}^{2}(f)=0$
 does not imply $X_{s}(f)=0$. Consider for example $f$ such that $\frac{\delta f}{\delta x}=\frac{\delta f}{\delta y}$.

  Assume $X(f)$ to $U^{\diamondsuit}$,
 then we have that if $f$ Hamiltonian, $d U^{\diamondsuit}(f)=0$ iff $\overline{\delta} (U^{\diamondsuit} f)=0$.
 Further, given $g$ Hamiltonian and $f \in \mathcal{D}_{L^{1}}$ with $\{ f,g \}=0$, we have 
 that $d U^{\diamondsuit}(f)=0$ iff $\overline{\delta} (U^{\diamondsuit} g)=0$. Thus, given 
 $(U^{\diamondsuit})^{\blacktriangle} f=0$ on a set of positive measure, we have that 
 $d U^{\diamondsuit} (f)= d (U^{\diamondsuit})^{\blacktriangle} (f)=0$.

 \newtheorem{seq}[adj]{Lemma Approximation with sequential movements and ac}
 \begin{seq}
 Given $\mathcal{G}$, we can associate $\Omega=\cup (x,y) \mid_{I}$ where $(U \gamma)_{I}=U_{j} \gamma_{I}$, 
 that is $U=U_{j}$ on a parameter interval $I$. 
 Given monotropy for a continuous movement $Uf$,
 there is within $\epsilon$ a sequential approximative movement. Given $U$ and $dU$ absolute continuous,
 we have that, given $\mid I \mid \neq 0$, we can approximate with a sequential movement. 
 We can determine the minimal $\mid I \mid$, such that movement
 does not change character on any smaller interval. For this, and smaller intervals, the movement 
 has constant dimension. Alternatively, let $C$ define a division of $\Omega$ and $I$. 
 Given a proper mapping $\psi : C \rightarrow \tilde{C}$ compact, assume $U_{j}^{\blacktriangle}=0$ 
 on $\tilde{C}$. Then
 $U^{\blacktriangle}_{j} \rightarrow 0$ is not necessarily regular limit, that is 
 $U_{j}^{\blacktriangle}$ is not absolute continuous.
 \end{seq}
 \subsection{Main result}
 
 Consider $H = F {}^{t} F$ and $U$ a movement with symmetric coefficients $\xi,\eta$.
 Further $(x,y) \rightarrow (x,z)$, where $z=y/x$. Given $y=y(\frac{1}{x})$ and 
 $y_{1}(x)=\frac{1}{x} y(\frac{1}{x})$, then $y_{1}$ can be seen as analytic outside 
 $\mid x \mid \geq R$, given $y$ analytic in $1/x$. Note that 
 ${}^{t} U y \rightarrow {}^{t} U^{\diamondsuit} y$ corresponds to $y(x) \rightarrow y(-1/\overline{x})$,
 when ${}^{t} U y$ is harmonic. In the following result, we assume the spiral $U$ is represented
 as a geometric mean $U \sim \sqrt{U_{1} U_{2}}$.

 \newtheorem{sumsofsq}[adj]{Proposition Sums of squares}
 \begin{sumsofsq}
 Consider the problem if the polar to $\Sigma X_{j}^{2}$ can be a spiral domain. 
 We have that $U^{*}$ spiral implies $U$ spiral. Given $H=F {}^{t} F$, we have that 
 ${}^{t} H=H$. When $H$ is considered over $\dot{B}$, that is $<H,\widehat{g}>$, where 
 $g \in \mathcal{D}_{L^{1}}$, we can assume ${}^{t} U \simeq U^{*} \simeq U^{\bot}$. Further,
 $(U F) {}^{t} (UF) \sim U F {}^{t} F U^{\bot}$. Thus, $<H,\widehat{g}> \sim <U H, U^{\bot} \widehat{g}>
 \sim <U^{2} H,\widehat{g}> \sim <U_{1} U_{2} H, \widehat{g}>$, where $U \sim \sqrt{U_{1} U_{2}}$,
 that is the spiral is well defined over $<H,\widehat{g}>$, given that $(U g) \in \mathcal{D}_{L^{1}}$.
 Given $<H,g>=0$ where $g \in \mathcal{D}_{L^{1}}$, we have when $H \in \mathcal{D}_{L^{1}}$ 
 that $g=0$ for $H$ hypoelliptic, that is the complement to the range is trivial. 
 For general $H$ as above, we have that $U^{\bot} g \bot U H$ or $U^{2} H \bot g$ 
 on a parameter interval. The condition $U^{\bot} g \bot U H$, means that $U^{\bot}g$ describes a polar.
 \end{sumsofsq}

 When the topology is such that 
 $\log H \in L^{1}$, we consider $U e^{\nu} \bot U e^{f^{2}} \sim e^{U \nu} \bot e^{U f^{2}}$
 and $\int_{\Sigma} e^{U \nu + U f^{2}} =0$ implies $U \nu + U f^{2}=0$ on $\Sigma$. 
 Assume $U$ harmonic and that $U f^{2} \sim (U f)^{2} \geq 0$. Given $g$ is defined such that 
 $g = \nu > 0$ on $\Sigma$ and given $U g$ has the same property, we see that $g=0$. 
 
\newtheorem{subideal}[adj]{Lemma Sums of squares is a proper sub ideal}
\begin{subideal}
 Given $\mathcal{G}(I) \subset (I)$, where $(I)=(I_{HE})$, assume we have $(I_{0}) \subset (I)$ with $F \in (I_{0})$, implies
 $F=\Sigma F_{j}^{2}$ and such that $F_{j}^{2} \in (I_{0})$ implies $U_{s}^{2} F_{j}^{2} \subset (I_{0})$. 
 Then $U_{s}^{2} F \in (I_{0})$ but $U_{s} \notin \mathcal{G}$.
 \end{subideal}

 Assume $H=F {}^{t} F \simeq (F,M)$, where $M=M(\nu,\vartheta)$. If $\Sigma$ is the set 
 where $H$ is hypoelliptic, assume $UH$ hypoelliptic on $\tilde{\Sigma}$, the continuation of $\Sigma$. 
 Assume $\tilde{C}$ the compact set, where $M(\nu,\vartheta)=M(\vartheta,\nu)$. 
 Assume $V$ such that $U H \simeq (F,VM)$ and $d \mu_{V}$, the measure such that 
 $V f = \int f d \mu_{V}$. Given $V$ absolute continuous,we thus have $Md \mu_{V}=0$ iff $V^{\blacktriangle}M=0$. 
 Transversality means $d \mu_{V} \bigoplus d \mu_{V^{\bot}}=0$ and 
 $\mbox{ supp } d \mu_{V} \cap \mbox{ supp } d \mu_{V^{\bot}}=\{ 0 \}$. Projectivity, means $V + V^{\blacktriangle}=I$.
 
 In the case when
 $V$ is not absolute continuous,we assume there is $U$ analytic on $M$, such that $\mid U - V \mid \leq \epsilon$ close to 
 $\mid \nu \mid=1$ and $d \mu_{U}$ is transversal.  Assume $\nu_{0}$ is a point on $\mid \nu \mid=1$, 
 where $\nu \rightarrow \nu_{0}$ regularly. 
 Assume $d M(\nu_{0},\vartheta)=0$. Given $d M$ analytic and $d M=0$
 on $\mid \nu \mid=1$, then $d M =0$ on $\mid \nu \mid \leq 1$ (or $\mid \vartheta \mid \leq 1$).
 Given $d M$ closed, we have $d M=0$ over $\nu \rightarrow \nu_{0}$ regularly. Assume
 $d U^{\bot} M=0$, then using that $\mathcal{G}$ very regular, there is a $V$ with
 ${}^{t} V M$ analytic close to $ \mid \nu \mid=1$ and $\mid U - V \mid \leq \epsilon$
 over $\mid \nu \mid=1$, that is $d V^{\bot} M =0$. Assume $E$ such that
 $E H - \delta_{0} \in C^{\infty}$, such that $d E \sim dM$, where we assume $d M \in L^{1}$
 and that $d M=0$ over a compact and symmetric set $K$. Then, assuming $d E$ very regular,
 we have that $d M$ can be represented regularly outside $\mid \nu \mid=1$. Note that when
 $K$ is represented using the \"ahnlich transform, when the sub level sets to $M$ are compact
 in $\nu$, they are not compact in $\vartheta$. Thus, when $d M$ has a regular representation,
 the (non-trivial) support is not symmetric.

 \section{Spectral projection}

 Given that $E_{\lambda}(v)=\int e_{\lambda}(x,y)v(y) d y$ a spectral projection, that is $I = \int d E_{\lambda}$,
 we have that $E_{\lambda}^{2}=E_{\lambda}$, 
 thus $\int e_{\lambda}(x,y) e_{\lambda}(y,z) d y=e_{\lambda}(x,z)$. Assume $X=X' \times X''$. 
 Given $e_{\lambda}(x',z')=\int_{f(\xi') < \lambda} e^{i (x'-z') \cdot \xi'} d \xi'$ and
 given $\int_{f < \lambda} d \xi''=1$, we see that $\int_{f < \lambda} e^{i(y'' - y'') \cdot \xi''} d \xi''=1$, 
 given that $\{ f < \lambda \}$ is compact. Thus, the $E_{\lambda}$ corresponding to a phe 
 operator, is projective on compact sub level sets. Alternatively we can consider the regularized
 spectral projector as in (\cite{Dahn15}).

 \newtheorem{dense}[adj]{Lemma Analytic subspace for graph of movement}
 \begin{dense}
 Assume $g \in \mathcal{D}_{L^{1}}$ and $\widehat{g} \in \dot{B}$, that is $U^{\bot} (=U^{*})  \in \mathcal{D}_{L^{1}} '$
 (algebraic base). Assume $\mathcal{G} \times \mathcal{G}^{\bot}$ very regular, 
 such that there is a subspace $H \times H$, where $(U,-U^{\bot})$
 preserves analyticity.  
 Assume $U^{\bot} \widehat{f}=0$, where $f,\widehat{f} \in H$ 
 and $H$ dense in 
 $\mathcal{D}_{L^{1}}$ (isolated singularities). Thus $U^{\bot} f=g$ implies $f=g$.
 \end{dense}
 
 Note that if $\Sigma^{*}$ is the set where $U^{*} \widehat{f} =0$, we have that
 $d U^{*}(\widehat{f}) =0$ on $\Sigma^{*}$ and since $U^{*} \neq U$, $d U(\widehat{f})=0$.
 Thus, if $U$ absolute continuous, $U^{\blacktriangle} \widehat{f}=0$. 
 Assume $<f,\widehat{g}> \sim <\widehat{f},g>$, then we have that 
 $< Uf, \widehat{g}> \sim <f, U^{\bot} \widehat{g}>$. Thus, we have that with respect to 
 $<, >$, that ${}^{t} U \simeq U^{\bot}$ and ${}^{tt} U \simeq U^{\bot \bot}$. The ideal 
 $(J) \ni h=\widehat{g} \in \mathcal{D}_{L^{1}}$ can be considered in a B-rum (cf. \cite{Nilsson72}) $\mathcal{B}_{\alpha}$).
 Given $f$ is considered in $H'$, we can consider $\widehat{f}$ in $Exp_{\parallel \cdot \parallel_{1}}$.
 
  Given $V^{\triangle}$ is a closed extension,
  with $N(V^{\triangle})=\{ 0 \}$ and sub harmonic on a Riemann surface $W \in \mathcal{O}_{G}$, 
  then we have that $V^{\triangle}=const. I$, note that $V^{\triangle}=const. I$ implies $V^{\blacktriangle}=0$.
  If for every $V^{\triangle}$ harmonic, we have $Flux (V^{\triangle})=0$ over the boundary, 
  we have that the domain limited by the boundary $\in \mathcal{O}_{G}$ and the boundary is removable. 
  for the movement (\cite{AhlforsSario60})

  Assume $U$ defined by $X(f)=\xi \frac{\delta f}{\delta x} + i \eta' \frac{\delta f}{\delta y}$.
  Assume $f \in \mathcal{D}_{L^{1}}$ and $\widehat{f} \in \dot{B}$. 
  Assume further, existence of $g \in H$, such that $U^{\bot} \widehat{g}$ analytic and 
  $U g=f$ (``first surface''). Given $U I$ a normal operator, we have that 
  $\parallel U^{\bot} g \parallel \sim \parallel U g \parallel$. Given Parseval
  we have $\parallel U^{\bot} \widehat{g} \parallel_{2} \sim \parallel U g \parallel_{2}$.

  \newtheorem{deficiency}[adj]{Definition Deficiency index} 
  \begin{deficiency}
  Define $\phi_{1} : U \rightarrow U^{\bot}$ and $\phi_{2} : (U,-U^{\bot}) \rightarrow (U,-U^{\bot})^{\bot}$.
  Thus $(U,-U^{\bot}) \in \mathcal{G} \times \mathcal{G}^{\bot}$. Given $\phi_{1}^{2} \sim id$
  that is $U^{\bot \bot} \sim U$, we have that $(U,-U^{\bot})^{\bot} \sim (U^{\bot},-U)$, that is $\phi_{2}(U,-U^{\bot})=(U^{\bot},-U)$.
  Thus $D(\phi_{2})$ (the domain) $=G_{\phi_{1}}$ (graph). Given $\phi_{1}$ symmetric, for instance
  $\phi_{1}(\overline{U})=\overline{\phi_{1}(U)}$, we have that 
  there is a continuation $H \rightarrow \mathcal{D}_{L^{1}}$ $\tilde{\phi}_{1}(U)=U^{\triangle}$.
  Further, $G_{\tilde{\phi_{2}}}=G_{\phi_{2}} \bigoplus \tilde{D}_{+} \bigoplus \tilde{D}_{-}$ and
  it is for a unique continuation, necessary to have defect index equal and zero.
  \end{deficiency}

 Concerning the extended plane, consider $\frac{d y}{d x} \rightarrow (X,Y) \rightarrow (x,y)$
 corresponding to dynamical systems, that is we assume $f$ Hamiltonian. When we consider 
 $F {}^{t} F \simeq (F,M)$, we have that as long as $\int d M$ is constant $=\mu$, 
 $\lambda \in \sigma(F)$ iff $\lambda \mu \in \sigma(F {}^{t} F)$.
 The eigen vectors in the extended plane, are given by $(\nu,\vartheta)$, such that 
 $ d M(\nu,\vartheta)$ bounded. Note that convexity with respect to $V$ implies in the 
 plane $\nu,\vartheta$, that the domain is on one side of a hyper plane.

 Consider $L^{2}=L^{2}_{ac} \bigoplus L^{2}_{p}$.
 Thus, when $\lambda \in \sigma(F {}^{t} F)$, we have that $L_{p}^{2} \neq \{ 0 \}$. 
 For hypoelliptic operators, we have that the localizer $F$, considered in $L^{2}$, 
 is very regular, that is $F^{\bot}$ has regularizing action.  
 Note that when $f^{N}$ hypoelliptic, we have that $F_{N}$ is very regular, that is $F_{N}^{\bot}$ has regularizing action.
 When $F_{N}^{\bot}$ has regularizing action,we do not have that 
 $F^{\bot}$ has regularizing action. However, we have that when $F_{N}$ very regular,
 then $F$ is very regular outside the kernel. 
 Thus, the localizer $F$ has a projective property outside the kernel (\cite{Dahn15}).
 
 \newtheorem{triweal}[adj]{Lemma Condition for trivial polar}
 \begin{triweal}
 When $F$ has kernel, we have that when $\phi \in \mbox{ ker }F$, $(\int M) F(\phi)=0$. 
 Note that when 
 $M={}^{t} M$ and the support for $M$ is one sided, we must have that $M$ has trivial support. 
 In the same manner if $d M$ is algebraic, we have that $\int_{\Omega} d M=0$ implies 
 $\lambda \Omega=0$.
 \end{triweal}
  
 Assume $V$ defined by $X(f)=0$, then there is a maximal domain 
 $\Omega$, such that $V f \mid_{\Omega}$ is in $L^{1}_{ac}$. Given $d U_{j}$ absolute continuous with respect to 
 $d U_{1}$ and $d U_{j}=\alpha_{j} d U_{1}$, we have that $\alpha_{j}$ has no zero's on a 
 maximal domain $\Omega$ (half space). Given $\Pi \alpha_{j} \neq 0$ we must have for all $\alpha_{j} \neq 0$,
 that is $\Sigma U_{j}$ absolute continuous.

 \newtheorem{spectral}[adj]{Lemma Spectral condition for projectivity}
 \begin{spectral}
 Given $f$ (formally) hypoelliptic, the spectral function is regularizing (\cite{Nilsson72}). 
 Assume $E_{\lambda}$ the projection corresponding to $f$ and $\tilde{E}_{\lambda}$ 
 corresponding to $Uf$.  Then $\tilde{E}_{\lambda}$ 
 is a projection operator if $U \rightarrow {}^{t} U$ projective. Conversely, when 
 $\tilde{E}_{\lambda}$ is projective, $U \rightarrow {}^{t} U$ preserves character.
 \end{spectral}
 
 When $U$ harmonic, 
 $(U E)^{2} \simeq U^{2} E^{2} \simeq U^{2} E \simeq U E$, 
 given that $U^{2} \simeq U$ (preserves character). However, when $U=U_{s} \simeq \sqrt{U_{1} U_{2}}$, we have
 $U^{2} \simeq U_{1}U_{2}$, that is $U^{2} \neq U$. The conclusion is that for $(UE)$ to be 
 a projection, we must assume $U \neq U_{s}$. Note that when $f$ and ${}^{t} U f$ are hypoelliptic and 
 $U$ reflexive and projective, we have that $UE$ is regularizing. Further, the condition $U I = I U$
 implies that $U \rightarrow {}^{t} U$ preserves character.
 
  Assume $F^{\bot}$ a projection operator (very regular) such that $F^{\bot}(\phi)=0$ implies $\phi =0$
  (modulo regularizing action). When $U$ is surjective, we have $F^{\bot} ({}^{t} U \phi)=0$ implies $\phi =0$.
  If we extend ${}^{t} U$ algebraically to $L^{1}$ using for instance 
  ${}^{t} U \rightarrow U \rightarrow U^{*} \rightarrow  U^{\triangle} \rightarrow  U^{\bot}$,
  the same conclusion holds for $U^{\bot}$. Thus, given that $F^{\bot}$ is completely
  determined by $U^{\bot} \phi$ ($({}^{t} U \phi)^{\bot} \simeq \{ 0 \}$), we have $UF + F^{\bot}U^{\bot}=I$,
  that is $U$ preserves projectivity. If we write $F^{\bot}$ as $(F^{\bot}, M)$, $M$ must
  have point support. Thus $U^{\bot} \phi = U \phi=0$ and when $U \rightarrow I$ regularly, we have
  $\phi=0$ ($=\phi(0)$).
 
 \section{Unique continuation property}
  
 Given $A f=X \frac{\delta f}{\delta x} + Y \frac{\delta f}{\delta y}$, there is a corresponding 
 $A^{*} \widehat{f}=\widehat{X} \frac{\delta \widehat{f}}{\delta x^{*}} + \widehat{Y} \frac{\delta \widehat{f}}{\delta y^{*}}$
 relating translation and rotation, such that $X (A f)=A (Xf)$ and thus in the same manner for 
 $X^{*}(A^{*} \widehat{f})=A^{*} (X^{*} \widehat{f})$

 \newtheorem{Lieortho}[adj]{Lemma The orthogonal relative the Lie algebra}
\begin{Lieortho}
Assume $\mathcal{G}=\mathcal{G}_{1} \bigoplus \mathcal{G}_{2}$, such that $U_{1} \sim U_{2}^{\bot}$. 
Then we have that $d \mathcal{G} = d \mathcal{G}_{1} \oplus d \mathcal{G}_{2}$. In particular 
when $I \in \mathcal{G}$, $d \gamma = d U_{1}(\gamma) \bigoplus d U_{2}(\gamma)$. Transversality 
means that $\Sigma_{j}= \{ d \mathcal{G}_{j}=0 \}$ has $\dim \Sigma_{1} \cap \Sigma_{2}=0$.  Given $\gamma$ 
polynomial (locally) and $\Sigma=\{ \gamma d \mathcal{G}=0 \}$, we have that $\dim \Sigma = \dim \{ d \mathcal{G}=0 \}$.
Note that it is necessary that $d U=0$ is given an orientation, for the normal to be locally 
algebraic.
\end{Lieortho}
 
 The continued group is derived in (\cite{Lie91} Chapter 13) through $\frac{d x}{d t}=\xi$,$\frac{d y}{d t}=\eta$
 and $\frac{d y'}{d t}=(\frac{\delta \eta}{\delta x} + (\frac{\delta \eta}{\delta y} - \frac{\delta \xi}{\delta x}) y' - \frac{\delta \xi}{\delta y} y^{'2})$.
 The infinitesimal transformation associated to the continued group is given by $X'(f)=\xi \frac{\delta f}{\delta x} + \eta \frac{\delta f}{\delta y} + \eta' \frac{\delta f}{\delta y'}$,
 where $\eta'=\frac{d y'}{d t}$.
  Concerning the definition of $\eta'$, assume $(\xi,\eta)$ related by the Cauchy-Riemann condition, 
  then we have that $\eta'=\frac{\delta \eta}{\delta x}(1 + y^{2'})$, that is if $y'$ is real 
  we have $\eta'=0$ iff $\frac{\delta \eta}{\delta x}=0$.
  Note that $y^{2'}$ algebraic does not imply that $y'$ is algebraic.
  For the continued equation, we thus have that $\tilde{X}(f) \equiv X(f)$, where 
  $\eta' \frac{\delta f}{\delta y'} \equiv 0$.

  We assume for $A^{\triangle} f$, that $\eta^{\triangle}/\xi^{\triangle} \sim \eta^{*} / \xi^{*}$, 
  where $(\xi^{\triangle},\eta^{\triangle})$ continuation according to Lie.  
  Given the continuation analytic, we have that $\eta^{' \triangle}=\frac{d \eta^{\triangle}}{d x^{*}} - \frac{d \xi^{\triangle}}{d y^{*}} y^{' 2}$.
  Note that $(\widehat{Y})_{x^{*}} - (\widehat{X})_{y^{*}} \neq 0$, thus if $(\widehat{Y}/\widehat{X})^{2}=1$,
  $\frac{d \eta^{\triangle}}{d x^{*}} - \frac{d \xi^{\triangle}}{d y^{*}} \neq 0$.
  
 Assume $\frac{\delta (V f)}{\delta x}=\xi \frac{\delta f}{\delta x}$ and $\frac{\delta (Vf)}{\delta y}=\eta \frac{\delta f}{\delta y}$,
then we have that $\Delta(Vf)= \xi \frac{\delta^{2} f}{\delta x^{2}} + \eta \frac{\delta^{2} f}{\delta y^{2}} +
\frac{\delta \xi}{\delta x} \frac{\delta f}{\delta x} + \frac{\delta \eta}{\delta y} \frac{\delta f}{\delta y}$.
Note that if $\xi/\eta \sim \frac{\delta f}{\delta x} / \frac{\delta f}{\delta y}$, we have that
$(\frac{\delta \eta}{\delta x} - \frac{\delta \xi}{\delta y})=\frac{\delta^{2} f}{\delta x \delta y} - \frac{\delta^{2} f}{\delta y \delta x}$

Consider $X(Af)=\xi \frac{d A}{d x} + \eta \frac{d A}{d y}$. On the edge 
$0=\frac{d^{2} F}{d T^{2}}=\frac{d^{2} F}{d x^{2}} (\frac{d x}{d T})^{2} + \frac{d F}{d x} \frac{d^{2} x}{d T^{2}} ´+ \frac{d^{2} F}{d y d x} \frac{d y}{d T} \frac{d x}{d T}
+ \frac{d^{2} F}{d y^{2}} (\frac{d y}{d T})^{2} + \frac{d^{2} F}{d x d y} (\frac{d x}{d T})(\frac{d y}{d T}) + \frac{d F}{d y} \frac{d^{2} y}{d T^{2}}$.
Given all second order derivatives to $F$ vanish, $\frac{d^{2} F}{d T^{2}}=\frac{d F}{d x} \frac{d^{2} x}{d T^{2}} + \frac{d F}{d y} \frac{d^{2} y}{d T^{2}}=0$.
Thus, $X(Af)=(\xi \frac{d X}{d x} + \eta \frac{d X}{d y}) \frac{d F}{d x} + (\xi \frac{d Y}{d x} + \eta \frac{d Y}{d y}) \frac{d F}{d y}$
 
 Concerning $X (Af)= \lambda A f$, assume $Y/X \sim -\eta / \xi$. Then, $X (Af)=\xi \frac{\delta A f}{\delta x} + \eta \frac{\delta A f}{\delta y}$
$\sim \frac{\delta f}{\delta x}(\xi X_{x} + \eta X_{y}) + \frac{\delta f}{\delta y}(\xi Y_{x} + \eta Y_{y})$ +
$\xi^{2} \frac{\delta f}{\delta x^{2}} - \eta^{2} \frac{\delta f}{\delta y^{2}} -
\xi \eta (\frac{ \delta f}{\delta x \delta y} - \frac{\delta f}{\delta y \delta x})$. Assume $X(X)=\lambda X$ and $X(Y)=\lambda Y$,
then the two first terms can be written $X(X) \frac{\delta f}{\delta x} + X(Y) \frac{\delta f}{\delta y} + \ldots$.
that is given the domain is such that $f_{xy}=f_{yx}$, it remains $\xi^{2} f_{xx} - \eta^{2} f_{yy}$.

 \newtheorem{polarpm}[adj]{Proposition On the projection method}
 \begin{polarpm}
 Assume the polar is defined as the set of $(d x,d y)$ such that $\mathcal{T}(\tilde{M},\tilde{W})=0$. 
 Consider $\widehat{(\tilde{M},\tilde{W})} \rightarrow (\tilde{M},\tilde{W})^{\triangle}(d \gamma)=$ 
 $(\tilde{M}_{2},\tilde{W}_{2})(d \gamma^{\triangle})$ continuous (\cite{Dahn13}). Thus. given the completion algebraic, 
 we have that $\widehat{(\tilde{M},\tilde{W})}(d \gamma)=0$ implies $ d \gamma^{\triangle}$ (closed) is in the polar.
\end{polarpm} 

\newtheorem{project}[adj]{Lemma On the projection method}
\begin{project}
Assume $M \sim -(X_{x} + Y_{y})$ and $W \sim Y_{x} - X_{y}$.
Given $(y')^{2} \sim 1$, we have that $\eta'=0$ implies $W + M y'=0$, that is $-W/M \sim \frac{dy}{dx}$ 
Thus $W/M \sim W_{2}/M_{2}$, corresponding to exact forms (\cite{Dahn13}).   
\end{project}

\cite{Treves79},\cite{Treves06},\cite{Hormander}
\bibliographystyle{amsplain}
\bibliography{sos}
\end{document}